\documentclass[11pt]{amsart}                    
\usepackage{amssymb}  

\input xy
\xyoption{all}
 
\theoremstyle{plain}
\newtheorem{theorem}{Theorem}[section]
\newtheorem{corollary}[theorem]{Corollary}
\newtheorem{lemma}[theorem]{Lemma}

\newtheorem*{conjecture}{Conjecture}

\newtheorem{fact}[theorem]{Fact}
\newtheorem{claim}[theorem]{Claim}

\theoremstyle{definition}
\newtheorem{definition}[theorem]{Definition}
\newtheorem{remark}[theorem]{Remark}

\newtheorem{convention}[theorem]{Convention}

\newtheorem{axiom}{Axiom}

\newcommand{\ran}{\operatorname{ran}}
\newcommand{\id}{\operatorname{id}}
\newcommand{\D}{\operatorname{Dep}}
\newcommand{\LS}{\operatorname{LS}}

\newcommand{\dom}{\operatorname{dom}}

\newcommand{\ftp}{\operatorname{tp}}

\newcommand{\gatp}{\operatorname{ga-tp}}
\newcommand{\tp}{\operatorname{ga-tp}}
\newcommand{\gaS}{\operatorname{ga-S}}

\renewcommand{\phi}{\varphi}

\newcommand{\Union}{\bigcup}

\newcommand{\initial}\lessdot
 \newcommand{\infinity}{\infty}
\newcommand{\K}{\operatorname{\mathcal{K}}}

\def\?{?\vadjust
{\vbox to 0pt{\vskip-7pt\hbox to 1.1\hsize{\hfill\huge ?!}}}}

\newbox\noforkbox \newdimen\forklinewidth
\setbox0\hbox{$\textstyle\smile$}
\setbox1\hbox to \wd0{\hfil\vrule width \forklinewidth depth-2pt
 height 10pt \hfil}
\setbox\noforkbox\hbox{\lower 2pt\box1\lower 2pt\box0\relax}
\def\unionstick{\mathop{\copy\noforkbox}\limits}

\def\nonfork_#1{\unionstick_{\textstyle #1}}

\setbox0\hbox{$\textstyle\smile$}
\setbox1\hbox to \wd0{\hfil{\sl /\/}\hfil}
\setbox2\hbox to \wd0{\hfil\vrule height 10pt depth -2pt width
               \forklinewidth\hfil}
\newbox\doesforkbox
\setbox\doesforkbox\hbox{\lower 2pt\box1 \lower 2pt\box2\lower2pt\box0\relax}
\def\nunionstick{\mathop{\copy\doesforkbox}\limits}

\def\fork_#1{\nunionstick_{\textstyle #1}}

\begin{document}

\title{Abstract decomposition theorem and
applications}

\author{Rami Grossberg}
\email[Rami Grossberg]{rami@cmu.edu}

 \address{Department of Mathematical Sciences\\
Carnegie Mellon University\\
Pittsburgh, PA 15213\\
United States}

\author{Olivier Lessmann}
\email[Olivier Lessmann]{lessmann@maths.ox.ac.uk}

\address{Mathematical Institute\\
Oxford University\\
Oxford OX1 3LB\\
United Kingdom}

\thanks{This is part of the second author's Ph.D. thesis, under the guidance
of the first author}
\date{February 9, 2004} 

\begin{abstract} 
In this paper, we prove a decomposition theorem
for abstract elementary classes $\K$ with the amalgamation property,
under the assumption that certain axioms regarding independence,
existence of some prime models, and regular types are satisfied.
This context encompasses the following:

\begin{enumerate}
\item
$\K$ is the class of models of 
an $\aleph_0$-stable first order theory.

\item
$\K$ is the class of $\mathbb F^a_{\aleph_0}$-saturated models of a
superstable first order theory.

\item
$\K$ is the class of models of 
an excellent Scott sentence $\psi\in L_{\omega_1,\omega}$.

\item
$\K$ the class of locally saturated models of a superstable good diagram $D$.

\item
$\K$ is the class of $(D,\aleph_0)$-homogeneous models of 
a totally transcendental good diagram $D$.
\end{enumerate}

We also prove the nonstructure part necessary to obtain a Main Gap theorem
for (5), which appears in the second author's Ph.D. thesis~\cite{Le}. 
The main gap in the contexts (1) and (2) are theorems of Shelah~\cite{Sh a};
(3) is by Grossberg and Hart~\cite{GrHa}; (4) by Hyttinen and 
Shelah~\cite{HySh:2}. 
\end{abstract}
\maketitle

\section*{Introduction}

This paper has two purposes.  The first is to present an abstract setting
lifting the essential features of {\em classifiable first order theories},
to settings which are not first order.
The second is to present, as an application, a new Main Gap theorem in the
context of homogeneous model theory.

In his celebrated paper~\cite{Sh 131}, Saharon Shelah proved the so-called 
\emph{Main Gap Theorem} for the class of $\aleph_\epsilon$-saturated
models of a complete first order theory $T$.
The result consists of showing that, if there are fewer than
the maximum number of nonisomorphic models of cardinality $\lambda> |T|$,
then the theory $T$ is superstable and satisfies NDOP,
every
$\aleph_\epsilon$-saturated
model has a decomposition in terms of an independent tree
of small models, and furthermore, the tree is well-founded. 
This implies that the number of nonisomorphic 
models in each cardinal is bounded by a slow growing function.
This exponential vs. slow growing dichotomy in the number of nonisomorphic
models is what is referred to as the main gap.

This main gap phenomenon is not limited to the first order case.
About six years after \cite{Sh 131}, Rami Grossberg and Bradd Hart \cite{GrHa}
proved the main gap for the class of models
of an excellent Scott sentence (see \cite{Sh:872} or Lessmann's
paper
in this volume  
for a definition of excellence).
The crucial property allowing a decomposition
is also NDOP.

Another non first order context with main gap phenomena is homogeneous
model theory.  
Homogeneous model theory was introduced by Shelah in \cite{Sh 3}.
It is quite general and includes first order logic, 
e.c. models, Banach space model theory, many generic constructions,
classes of models with amalgamation over sets (infinitary, $L^n$),
as well as concrete cases like expansions of Hilbert spaces,
and free groups. 
We have good notions of 
omega-stability/total transcendence~\cite{Le1}, 
superstability (\cite{HySh:2}, \cite{HyLe}), 
stability (\cite{Sh:54}, \cite{Gr1}, \cite{Gr2}, \cite{GrLe}, \cite{HySh:1}), 
and simplicity~\cite{BuLe}.
For an exposition on homogeneous model theory, see \cite{Le:4}.
Hyttinen and Shelah proved a Main Gap theorem
for locally saturated models of a superstable diagram~\cite{HySh:2}~(4);
NDOP is also the dividing line.  Our result for the class of 
$(D,\aleph_0)$-homogeneous models of a totally transcendental diagram~(5)
also uses NDOP.
These results were proved independently and are incomparable 
approximations to the main gap conjecture for 
the class of 
$\aleph_\epsilon$-saturated models (it was proved under simplicity
by Hyttinen and Lessmann).

The Main Gap proofs have two components:  On the one hand we have
a structure result (each model is decomposed into an independent
tree of small models), and on the other hand we have a non-structure
result (undesirable features, for example if the tree is not well-founded,
produce the maximal number of models).
In this paper, we present a framework in which we prove the structure part.
This framework is defined in the context of Shelah's Abstract Elementary
Classes~\cite{Sh:88} 
with the Amalgamation Property and the Joint Embedding Property.
See \cite{Gr:3} for an exposition on this subject.
We postulate the existence of an independence relation which is well-behaved
over models.
We also postulate the existence of a special kind of prime models,
called primary, which exist over certain sets and behave well
with respect to the independence relation (dominance).
Finally, we posit the existence of certain types, called regular,
which must be dense, and we ask for good behaviour (essentially capturing
basic orthogonality calculus).
This framework is general enough to include
in the same proof \cite{Sh 131}, \cite{GrHa}, \cite{HySh:2},
and the case of $(D, \mu)$-homogeneous
models of a totally transcendental diagram $D$
introduced in \cite{Le1}.

In addition to providing a useful result allowing one to bypass all
the technical aspects of a structure proof,
we try to isolate and understand 
the essential features that ensure a good structural
theory.  
The motivation is similar to Shelah's
{\em good frames} (see \cite{Sh 600}, \cite{Sh:705}),
in his work around categoricity for abstract elementary classes.  
Here, we
think that  a decomposition theorem under NDOP is a good 
indication that a good structure
theory is possible.  
The next step is to see if geometric model theory can be developed in 
this context.  There are also other good indications that this is possible;
a recent result of the second author with 
Hyttinen and Shelah \cite{HLS} generalises
one of the basic results of first order geometric model theory obtained
for (1) -- (2) to (3) -- (5); we believe that generalisations to 
the context we isolate is possible.

Our result is a modest 
step towards the following conjecture of Shelah (late 1990s):
\begin{conjecture}[Shelah]  
Let $\K$ be an AEC.  Denote by $\delta$ the ordinal
$(2^{\LS(\K)})^+$.  If $\K$ has at least one model, but fewer 
than the maximal number of models
in some cardinal $\lambda>\beth_\delta$, then the number 
of nonisomorphic models of size $\aleph_\alpha$ is bounded by 
$\beth_\delta(\aleph_0+|\alpha|)$ for each ordinal $\alpha$.
\end{conjecture}

The paper is organised as follows:
In Section 1, we introduce an axiomatic framework for
abstract elementary classes with AP and JEP.
We have axioms postulating the existence of a good independence
relation capturing the essential features of the superstable case,
in a spirit similar to Baldwin~\cite{Ba}.
We have axioms on primary models, their existence over certain sets,
their uniqueness, and their behavior with respect to the independence
relation.
We also have axioms regarding the existence of regular types
and how they connect with independence and primary models.
The axioms are numbered separately for easy reference.
We prove a decomposition theorem in this axiomatic framework
under NDOP (Theorem~\ref{dec}).
The key difference with Shelah's abstract treatment of his main 
gap theorems~\cite{Sh b} is that his relies on compactness, whereas
ours does not.
We also describe how (1) -- (4) fall within this framework.

In Section 2, we present the necessary orthogonality calculus
to show that the class of $(D,\aleph_0)$-homogeneous models
of a totally transcendental diagram $D$ satisfies
the axioms of Section 1. 
This implies that under NDOP, every $(D, \aleph_0)$-homogeneous
model is prime and minimal over an independent tree of small
models.
We also prove several additional lemmas that
will allow us to complete the main gap for this class.

In Section 3, we introduce DOP (the negation of NDOP)
for the class of $(D,\aleph_0)$-homogeneous models
of a totally transcendental diagram $D$.
We show that DOP implies the existence of many nonisomorphic models
(Theorem \ref{dopmany}).
For nonstructure results using DOP (the failure of NDOP),
the axiomatization needs several levels of saturation
(or homogeneity, or fullness). 
We give a proof of the nonstructure parts of the theorem in the context 
of Chapter IV (of \cite{Sh b}).
This gives the main gap for the class $\mathcal{K}$ of $(D, \mu)$-homogeneous
models of a totally transcendental diagram $D$ (for
any infinite $\mu$).
Note that, since finite diagrams generalise the first order
case, it is easy to see that the failure of a finite diagram to be totally
transcendental does not imply the existence of many
models. 
All the basic tools in place, we can also show,
using the methods of \cite{Sh b} or \cite{Ha} 
that $\lambda \mapsto I(\lambda, \mathcal{K})$
is weakly monotonic (Morley's Conjecture) for sufficiently
large $\lambda$.

In Section 4, we introduce depth 
for the class of $(D,\aleph_0)$-homogeneous models.
We prove that if a class is deep then it has many nonisomorphic
models (Theorem \ref{deepmany}). 
Finally, we derive the main gap (Theorem \ref{maingap})
for this class.
Using the same methods, we can also derive the main gap
for the class of $(D,\mu)$-homogeneous models
of a totally transcendental diagram $D$.

We would like to thank John T. Baldwin for several useful comments.

\section{The axiomatic framework and decomposition theorem} 

We fix
$(\mathcal{K},
\prec)$ an {\em abstract elementary class}~\cite{Sh:88} (AEC for short), 
{\em i.e.}
we assume that $(\K, \prec)$ satisfies the following axioms:

\begin{definition}[Abstract Elementary Class]\label{a:aec}
$\K$ is a class of models in the same similarity type $L$.
The relation $\prec$ is a partial order on $\K$.
\begin{enumerate}
\item 
$\K$ is closed under isomorphism;
\item
If $M, N \in \K$ and $M \prec N$ then $M \subseteq N$,
{\em i.e.} $M$ is a submodel of $N$;
\item
If $M, N, M^* \in \K$ with $M \subseteq N$ and $M, N \prec M^*$ then
$M \prec N$;
\item
There is
a cardinal $\LS(\K)$ such that for all $M \in \K$ and $A \subseteq M$
there is $N \prec M$ containing $A$ of size at most $|A|+ \LS(\K)$.
\item
$\K$ is closed under Tarski-Vaught chains:
Let $(M_i : i < \lambda)$ be a $\prec$-increasing and continuous
chain of models of $\K$.  Then $\bigcup_{i < \lambda} M_i \in \K$.
Also $M_0 \prec\bigcup_{i < \lambda} M_i$ and further, 
if $M_i \prec N \in \K$ for each $i < \lambda$, then
$\bigcup_{i < \lambda} M_i \prec N$.
\end{enumerate}
\end{definition}

It is not difficult to see that if $I$ is a directed partially ordered
set and $(M_s : s \in I)$ is such that $M_s \in \K$ and $M_s \prec M_t$
for $s < t$ in $I$, then $\bigcup_{s \in I} M_s \in \K$.

We now define naturally {\em $\K$-embedding} to be those embeddings 
preserving $\prec$.
We will work under the additional hypothesis that $\K$ has
the {\em Joint Embedding Property} (JEP) and 
the {\em Amalgamation Property} (AP).

\begin{axiom}[{\bf Joint Embedding Property}]\label{a:jep}
Let $M_0, M_1 \in \K$.
Then there is $M^* \in \K$ and $\K$-embeddings 
$f_\ell : M_\ell \rightarrow M^*$ for $\ell = 1,2$.
\end{axiom}

\begin{axiom}[{\bf Amalgamation Property}]\label{a:ap}
For $M_\ell \in \K$ ($\ell=0,1,2$) and 
$\K$-embeddings $f_\ell : M_0 \rightarrow M_\ell$,
for $\ell = 1,2$, there exist $M^* \in \K$ and $\K$-embeddings 
$g_\ell : M_\ell \rightarrow M^*$, for $\ell = 1,2$, 
such that $g_1 \circ f_1 = g_2 \circ f_2$.
In other words, the following diagram commutes:
\[
\xymatrix{M_1 \ar[r]^{g_1} & M^{*} \\
M_0 \ar[u]^{f_1} \ar[r]_{f_2} & M_2 \ar[u]_{g_2}
}
\]
\end{axiom}

\begin{remark}  
Notice that when $\K$ is an AEC then the amalgamation property is equivalent
to:  For every $M_\ell\in \K$ (for $\ell=0,1,2$) such that 
$M_0\prec M_\ell$ (for $\ell=1,2$) there are
$M^* \in \K$ with $M_2 \prec M^*$  and a $\K$-embedding 
$g:M_1\rightarrow M^*$ such that the following diagram commutes:
\[
\xymatrix{M_1 \ar[r]^{g} & M^{*} \\
M_0 \ar[u]^{\id} \ar[r]_{\id} & M_2 \ar[u]_{\id}
} 
\]
\end{remark}

Recall the next definition:

\begin{definition}
$M \in \K$ is {\em $\lambda$-model homogeneous} if whenever $N_0 \prec M$
and $N_1 \in \K$ with $N_0 \prec N_1$ and $N_1$ has size less than $\lambda$,
then exists a $\K$-embedding $f: N_1 \rightarrow M$ which is the identity
on $N$.
\end{definition}
Assuming that $\K$ has AP and JEP, it is possible
to construct {\em $\lambda$-model homogeneous} models 
for arbitrarily large $\lambda$.
Notice that since we are not assuming the existence of large models in $\K$,
a $\lambda$-model homogeneous model $M$ may be small, even though $\lambda$
is big.
If $M$ is $\lambda$-model homogeneous, 
then any $N \in \K$ of size less than $\lambda$ $\K$-embeds into $M$.
We can then use model-homogeneous models as monster models
(more on this later).
We make the following convenient convention.

\begin{convention}
We fix a $\bar{\kappa}$-model homogeneous model $\mathfrak{C} \in \K$,
for a suitably large cardinal $\bar{\kappa}$.
We will work inside $\mathfrak{C}$; all sets and models are assumed
to be inside $\mathfrak{C}$ of size less than $\bar{\kappa}$.
\end{convention}

We now postulate the existence of an \emph{independence relation} on subsets
of $\mathfrak{C}$, {\em i.e.} 
a relation on triples of sets $A, B$, and $C$, 
\[
A \nonfork_B C,
\]
satisfying the some axioms.
These are similar to the standard first order axioms for non-forking
in the context of superstability; {\em i.e.} Local Character is 
with respect to a finite set.
They are weaker in one respect:  we deal mainly with models.
This appears in the phrasing of Symmetry, Transitivity, and Local Character.
We do not assume Extension; this will be done addressed later when
we deal with stationarity.

\begin{axiom}[{\bf Independence}] \label{a:indep} 
Let $A, B, C$ and $D$ be sets.
Let $M$ be a model.
\begin{enumerate}

\item 
({\bf Definition}) $A \nonfork_B C$ if and only if $A \nonfork_B B \cup C$;

\item 
({\bf Triviality})  Let $M$ be a model. 
Then
$A \fork_M A$;

\item 
({\bf Finite Character}) $A \nonfork_B C$ if and only if $A' \nonfork_{B} C'$,
for all finite $A' \subseteq A$,
$C' \subseteq C$;

\item 
({\bf Monotonicity}) 
If $A \nonfork_B C$ and 
$B \subseteq B_1 \subseteq M_1$ and $C' \subseteq C$, 
then $A \nonfork_{B_1} C'$;

\item
({\bf Local Character})
Let $(M_i : i < \omega)$ is an $\prec$-increasing sequence of models
and $M = \bigcup_{i < \omega} M_i$.
Then, for every $\bar{a}$ there is $i < \omega$ such that
$\bar{a} \nonfork_{M_i} M$;

\item 
({\bf Transitivity}) 
If $M_0 \subseteq M_1 \subseteq C$, then 
$A \nonfork_{M_1} C$ and $A \nonfork_{M_0} M_1$ if and only if 
$A \nonfork_{M_0} C$;

\item 
({\bf Symmetry over models}) 
$A \nonfork_M C$ if and only if $C \nonfork_M A$.

\item 
({\bf Invariance})
Let $f$ be a $\K$-embedding 
with $A \cup B \cup C \subseteq \dom(f)$.
Then $A \nonfork_B C$ if and only if
$f(A) \nonfork_{f(B)} f(C)$.

\end{enumerate}
\end{axiom}

We are now concerned with {\em prime} models.

\begin{definition}
We say that a model $M \in \K$ is \emph{prime over $A$}, if
for every $N \in \K$ containing $A$, there exists a $\K$-embedding
$f \colon M \rightarrow N$, which is the identity on $A$.
\end{definition} 

We work with a special kind of prime models,
called {\em primary} models.
We isolate the main property of {\em bona fide}
primary models that we are going to use,
namely that 
any two primary models over the same set are isomorphic (but
we do not assume that all prime models have this property):

\begin{axiom}[{\bf Uniqueness of primary models}]\label{a:primary} 
Let $M \in \K$ be primary over $A$.
Then $M$ is prime over $A$.
Moreover, if $M' \in \K$ is another primary model over $A$,
then $M$ and $M'$ are isomorphic over $A$.
\end{axiom}

First, we need to define the notion of independent system.
We will say that a set of finite sequences $I$ is
a \emph{tree} if it is closed under initial segment.
We will use the notation $\eta \prec \nu$ to mean
that $\eta$ is an initial segment of $\nu$.

\begin{definition}
Let $I$ be a tree, 
we say that $\langle M_\eta \mid \eta \in I \rangle \subseteq M_1$ is 
a \emph{system} if $M_\eta \in \K$ for each $\eta \in I$ and
$M_\eta \subseteq M_\nu$ when $\eta \prec \nu \in I$.
\end{definition}

The concept in the next definition is called \emph{system in stable
amalgamation} by Shelah (see \cite{Sh:872} and \cite{Sh b}). 

\begin{definition} 
We say that $\langle M_\eta \mid \eta \in I \rangle$ is 
an \emph{independent system} if
it is a system satisfying in addition:
\[
M_\eta \nonfork_{M_{\eta^-}} \Union_{\eta \not \prec \nu } M_\nu,
\quad
\text{for every $\eta \in I$}.
\]
Where $\eta^-$ is the predecessor of $\eta$, i.e. $\eta^-:=\eta\restriction
(\ell(\eta)-1)$.
\end{definition}

\begin{axiom} [{\bf Existence of primary models}] \label{c4:a:prime} \hfill

\begin{enumerate} 
\item 
Let $M \in \K$.
There is a primary model $M' \prec M$ over the empty set;

\item 
If $\bar{a} \in N \setminus M$ (where $\bar{a}$ is finite)
then there is a primary model $M' \prec  N$ over $M \cup \bar{a}$;

\item 
If $\langle M_\eta \mid \eta \in I \rangle \subseteq N$ is an 
independent system, then
there exists a primary model $M' \prec N$ over 
$\Union_{\eta \in I} M_{\eta}$.
\end{enumerate}
\end{axiom}

The next axiom is one of the key properties of primary models.

\begin{axiom} [{\bf Dominance}] \label{c4:a:dominance} 
Suppose that $A \nonfork_M C$ and $M(C)$ is primary over $M \cup C$.
Then $A \nonfork_M M(C)$.
\end{axiom}

This implies the Concatenation property of independence.
We prove it in detail, but will be terser in the future.
\begin{lemma}[Concatenation]
If $A \nonfork_M BC$ and $C \nonfork_M B$ then $C \nonfork_M BA$. 
\end{lemma}
\begin{proof}
By Finite Character, we may assume that $A, B$, and $C$ are finite.
Let $M' \prec M''$ with $M'$ primary over $M \cup B$ and $M''$ primary
over $M \cup BC$.
By Dominance, we have $A \nonfork_M M''$, so $A \nonfork_M M'C$ by 
Monotonicity, $A \nonfork_{M'} C$ by Definition and Monotonicity, 
so $C \nonfork_{M'} A$
by Symmetry.  
By Dominance, we also have $C \nonfork_M M'$.
Thus $C \nonfork_M M'A$ by Transitivity, so $C \nonfork_M BA$ by Monotonicity.
\end{proof}

\begin{lemma}\label{l:up}
Let $M \prec N \in \K$ and assume that $ab \nonfork_M N$.
Then $a \nonfork_M b$ if and only if $a \nonfork_N b$.
\end{lemma}
\begin{proof}
By Monotonicity, $a \nonfork_M N$, so if $a \nonfork_N b$, then
also $a \nonfork_M b$ by Transitivity.

For the converse, choose a primary model $M'$ over $M \cup b$.
we claim that $a \nonfork_{M'} Nb$.
To see this, notice that $N \nonfork_M ab$, so by Dominance
$N \nonfork_M M''$,
where $M''$ is a primary model over $M \cup ab$, which 
implies $N \nonfork_M M' a$ and $N \nonfork_{M'} a$ by Monotonicity,
so the claim follows by Symmetry.
Now assume that $a \nonfork_M b$.
By Dominance, we have $a \nonfork_M M'$, so $a \nonfork_N b$
by Transitivity using the claim.
\end{proof}

\begin{definition} 
We say that $\{B_i \mid i< \alpha\}$ is \emph{independent over $M$} if 
\[
B_i  \nonfork_M \Union \{ B_j \mid j\not = i, j<\alpha \},
\]
for every $i < \alpha$.
\end{definition}

We now prove the usual result on  
independent families.

\begin{lemma} \label{c4:l:indepfamily}
Let $\{ B_i \mid i< \alpha \}$ be a family of sets and
assume that
\begin{equation} \tag{*}
B_{i+1} \nonfork_M \bigcup \{ B_j \mid j < i \},
\quad
\text{for every $i < \alpha$}.
\end{equation}
Then $\{ B_i \mid i< \alpha \}$
is independent over $M$.
\end{lemma}

\begin{proof}
By finite character of independence, it is enough to prove this statement
for $\alpha$ finite.
We do this by induction on the integer $\alpha$.

For $\alpha=1$, (*) implies that $B_1 \nonfork_M B_0$,
so by symmetry over models we have $B_0 \nonfork_M B_1$,
which shows that $\{ B_0, B_1 \}$ is independent over
$M$.

Assume by induction that the statement is true for $\alpha < \omega$.
Let $i \leq \alpha+1$ be given. 
We must show that
\begin{equation} \tag{**}
B_i \nonfork_M \bigcup \{ B_j \mid j \leq \alpha+1, j \not = i \}.
\end{equation}
If $i = \alpha+1$, then this is (*), 
so we may assume that $i \not = \alpha+1$ and therefore
(**) can be rewritten as
\[
B_i \nonfork_M \bigcup \{ B_j \mid j \leq \alpha, j \not = i \} 
\cup B_{\alpha+1}.
\]
Notice that by induction hypothesis
\[
B_i \nonfork_M \bigcup \{ B_j \mid j \leq \alpha, j \not = i \}.
\]
And further by (*) we have
\[
B_{\alpha +1} \nonfork_M \bigcup \{ B_j \mid j \leq \alpha, j \not = i \}
\cup B_i.
\]
Hence, the result follows by Concatenation.
\end{proof}

Under our axioms, independent systems are quite independent. 
For $J$ a subtree of $I$, denote $M_J = \bigcup \{ M_\eta : \eta \in J \}$
(note that $M_J$ is not necessarily a model).
The following lemma is a version of Shelah's generalised Symmetry
Lemma.  It also appears in a similar form in Makkai~\cite{Mak}.

\begin{lemma}
Let $\langle M_\eta \mid \eta \in I \rangle$
be an independent system.
Then, for any $I_1, I_2$ subtrees of $I$, we have:
\begin{equation} \tag{*}
M_{I_1} \nonfork_{M_{I_1 \cap I_2}} M_{I_2}
\end{equation}
\end{lemma}
\begin{proof}
By the finite character of independence, 
it is enough to prove (*) for finite trees $I$.
We prove this by induction on $|I_1 \cup I_2|$.

First, if $I_2 \subseteq I_1$, then it is obvious.
Thus, assume that there is $\eta \in I_2 \setminus  I_1$,
and choose $\eta$ of maximal length.
Let $J_2 := I_2 \setminus  \{ \eta \}$.
Notice that by choice of $\eta$, we have 
$M_{I_1 \cap J_2} = M_{I_1 \cap I_2}$.
By induction hypothesis, we have that
\begin{equation} \tag{*}
M_{I_1} \nonfork_{M_{I_1 \cap I_2}} M_{J_2}.
\end{equation}
Since $M_{\eta^-} \subseteq M_{J_2}$, 
by monotonicity
(*) implies that
\begin{equation} \tag{**}
M_{I_1} \nonfork_{M_{\eta^-}} M_{J_2}.
\end{equation}
By definition of independent system
and monotonicity
we have
\begin{equation} \tag{***}
M_\eta \nonfork_{M_{\eta^-}} M_{I_1} \cup M_{J_2}.
\end{equation}
Therefore, by concatenation
applied to (**) and (***),
we can conclude that
\begin{equation}\tag{$\dagger$}
M_{I_1} \nonfork_{M_{\eta^-}} M_{I_2}.
\end{equation}
Now, using (*) and monotonicity 
we have 
\begin{equation} \tag{$\ddagger$}
M_{I_1} \nonfork_{M_{I_1\cap I_2}} M_{\eta^-}.
\end{equation}
Thus, the transitivity property
applied to ($\dagger$) and ($\ddagger$), implies
that
\[
M_{I_1} \nonfork_{M_{I_1 \cap I_2}} M_{I_2}.
\]
This finishes the proof.
\end{proof}

in the context of an abstract elementary class with amalgamation,
Shelah introduced a natural notion of {\em types over models}.
The material we are about to cover can be found in more details in \cite{Gr:3}.
Consider the following equivalence relation $\sim$
on triples $(a, M, N)$, where $M, N \in \K$, $M \prec N$
and $a \in N$.
We say that 
\[
(a_1, M, N_1) \sim (a_2, M, N_2)
\]
if  there is $N^*\in \K$ and $\K$-embeddings
$f_\ell : N_\ell \rightarrow N^*$ which are the identity on $M$
and such that $f_1(a_1) = f_2(a_2)$.
The picture is
\[
\xymatrix{\ar @{} [dr] N_1
\ar[r]_{f_1}  &N^* \\
M \ar[u]^{\id} \ar[r]_{\id} 
& N_2 \ar[u]_{f_2} 
}
\]

The amalgamation property ensures that $\sim$ is an equivalence
relation.
The equivalence class $(a,M,N)/\sim$ is written $\gatp(a/M,N)$;
it is {\em the galois type of $a$ over $M$ in $N$}.

\begin{definition}
Let $M \in \K$.  Then
\[
\gaS(M) = \{ \gatp(a/M,N) : 
\text{For some $a \in N$, $M \prec N \in \K$} \}.
\]
\end{definition}
As is common in first order model theory, we use the letters
$p, q$, and $r$ for types.
We will say that $N'$ {\em realizes} $\gatp(a/M,N)$
if $M \prec N'$ and there is $b \in N'$ such
that $\gatp(b/M,N')=\gatp(a/M, N)$.
We continue to write $b \models p$, if $a$ realizes $p$.
Similarly, we can define $p \restriction M'$ for $M' \prec M$,
and $p \subseteq q$.
The amalgamation property guarantees that the types are well-behaved;
for example an increasing and continuous union of types is a type.

We have the following striking correspondence between 
{\em $\lambda$-saturation} and $\lambda$-model homogeneity.
Recall first:
\begin{definition}
Let $M \in \K$.
Then $M$ is {\em $\lambda$-galois saturated} if $M$ realizes
each $p \in \gaS(N)$
for $N \prec M$ of size less than $\lambda$.
\end{definition}

The next fact is due to Shelah (see \cite{Sh 576}, also \cite{Gr:3} for a proof).
\begin{fact}[Shelah] 
Let $\K$ be an abstract elementary class with AP and JEP.
For any $\lambda>\LS(\K)$ we have that  
$M$ is $\lambda$-model homogeneous if and only if $M$
is 
$\lambda$-galois saturated.
\end{fact}

This justifies further our use of $\mathfrak{C}$ as a monster
model: All relevant types are realized in $\mathfrak{C}$.
From now on, since we may only consider types of the form 
$\gatp(a/M,\mathfrak{C})$, we 
will omit $\mathfrak{C}$.

The invariance of the independence relation makes it natural to extend
the independence relation to types.

\begin{definition} Let $M \in \K$.
\begin{enumerate}
\item 
We say that $p \in \gaS(M)$ is \emph{free over $N  \prec M$} 
if for every $\bar{a} \in \mathfrak{C}$ realizing $p$, 
we have $\bar{a} \nonfork_N M$;

\item 
We say that $p \in \gaS(M)$ is \emph{stationary} if 
for every $N \in \K$ containing $M$, there is a unique
extension $p_N \in \gaS(N)$ of $p$ such that $p_N$ is
free over $M$.

\item
We say that the stationary type $p \in \gaS(M)$ is \emph{based
on $N$} if $p$ is free over $N$.

\end{enumerate}
\end{definition}
It is clear that `for every' is equivalent to `for some' in (1).
Notice that in this context, the existence of a free extension to
a stationary type (2), not just its uniqueness, is quite important.

\begin{axiom}[{\bf Existence of Stationary types}] \label{c4:a:stat}
Let $M \in \mathcal{K}$.
Then any non-algebraic  $p \in \gaS(M)$ is stationary.
\end{axiom}

The next lemma follows from the definition, Local Character, and
Transitivity.

\begin{lemma}
Let $p \in \gaS(M)$ and let $(M_i : i < \lambda)$ be an $\prec$-increasing
and continuous chain of models such that $\bigcup_{i < \lambda} M_i = M$.
Then there is $i < \lambda$ such that $p$ is based on $M_i$.
In particular, there is always $N  \prec M$ of size $\LS(\K)$ such
that $p$ is based on $N$.
\end{lemma}

We now introduce a strong independence between 
stationary types: {\em orthogonality}.

\begin{definition} 
Let $p \in \gaS(M)$ and $q \in \gaS(N)$.
We say that $p$ is \emph{orthogonal} to $q$, written $p \perp q$, if
for every $M_1 \in \mathcal{K}$ containing $M \cup N$
and for every $a \models p_{M_1}$ and $b \models q_{M_1}$,
we have $a \nonfork_M b$.
\end{definition}

By symmetry of independence, $p \perp q$ if and only
if $q \perp p$.
Also, if $p \in \gaS(M)$, $q \in \gaS(N)$ and 
with $M \prec N$, then by definition
$p \perp q$ if and only if $p_N \perp q$.
In fact, more is true:

\begin{lemma} [Parallelism]
Let $M \prec N$ and $p, q \in \gaS(M)$.
Then $p \perp q$ if and only if $p_N \perp q_N$.
\end{lemma}
\begin{proof}
We have already shown the left to right direction, so 
suppose that $p_N \perp q_N$ and let $M_1 \in \K$
with $M \prec M_1$, and suppose, for a contradiction
that $a \models p_{M_1}$, $b \models q_{M_1}$, but
$a \fork_{M_1} b$.
Let $N_1 \in \K$ containing $M_1 \cup N$.
By stationarity, there exists $a'b' \models \tp(ab/N_1)$
such that $a'b' \nonfork_{M_1} N_1$.
Notice that $a' \fork_{M_1} b'$, so $a' \fork_{N_1} b'$ by
Lemma~\ref{l:up}.
But $a' \models p_N$ and $b' \models q_N$, contradicting
the fact that $p_N \perp q_N$.
\end{proof}

We now expand this definition to orthogonality against models.

\begin{definition} Let $M, N, M_\ell \in \K$ for $\ell= 0,1,2$.
\begin{enumerate}
\item
Let $p \in \gaS(N)$. 
We say
that $p$ is \emph{orthogonal to $M$},
written $p \perp M$, if $p$ is orthogonal
to each $q\in \gaS(M)$.

\item 
If $M_0 \prec M_1,  M_2$, we write
that $M_1 / M_0 \perp M_2$ if and only if 
$p \perp M_2$,
for every $p \in \gaS(M_0)$ realized in $M_1$.
\end{enumerate}
\end{definition}

We now concentrate on a special kind of types: {\em regular types}.

\begin{definition} 
A stationary type $p \in \gaS(M)$ is called \emph{regular} if
for any $N \prec M$ with $p$ based on $N$ 
and for any $M_1 \in \K$ containing $M$ and 
$q \in \gaS(M_1)$ extending $p \restriction N$, 
either $q = p_{M_1}$ or
$q \perp p$.
\end{definition}

\begin{lemma} \label{c4:l:regular}
Let $M \subseteq M_1$.
If $p \in \gaS(M)$ is regular,
then $p_{M_1} \in \gaS(M_1)$ is regular.
\end{lemma}
\begin{proof}
Let $p \in \gaS(M)$ be regular. 
Let $N \prec M$ with $p$ based on $N$.
Then, $p_{M_1}$ is stationary based on $N$.
Let $q \in \gaS(N_1)$ extend $p \restriction N$.
Then either $q = p_{N_1} = (p_{M_1})_{N_1}$ or $q \perp p$.
Hence by definition of $\perp$ we have $q \perp p_{M_1}$.
This shows that $q$ is regular.
\end{proof}

The next axioms guarantee that it is enough to focus on
regular types.
\begin{axiom} [{\bf Existence of Regular types}] \label{c4:a:existence}
If $M \subseteq N$ and $M \not = N$, 
then there exists a regular type $p \in \gaS(M)$
realized in $N \setminus M$.
\end{axiom}

\begin{axiom}[{\bf Perp I}] \label{c4:a:perpM} 
Let $M, N \in \mathcal{K}$ such that $M \prec N$. 
Let $p \in \gaS(N)$ be regular.
Then $p \perp M$ if and only if  
$p \perp q$, for every regular type
$q \in \gaS(M)$.
\end{axiom}

This is to establish connections with the dependence
relation and orthogonality.

\begin{axiom}[{\bf Equivalence}]\label{c4:a:perp}
Let $M \in \mathcal{K}$  
and let $p, q \in \gaS(M)$ be regular
and let $\bar{b} \not \in M$ realize $p$. 
Then $q$ is realized in $M(\bar{b}) \setminus M$ if and only
if $p \not \perp q$.
\end{axiom}

Note that by Equivalence, the relation $\not \perp$ among regular
types (over the same base set) is an equivalence relation.

\begin{lemma} \label{c4:l:equiv} 
Let $M_0 \subseteq M \subseteq M' \subseteq N$.
Let $p \in \gaS(M')$ be regular realized in $N \setminus M'$
and $q \in \gaS(M)$ such that $p \not \perp q$.
Let $r \in \gaS(M_0)$ be regular. 
If $p \perp r$ then $q \perp r$.
\end{lemma}
\begin{proof}
By Lemma \ref{c4:l:regular}, the types $r_{M'}$
and $q_{M'}$ are regular.
By definition, $p \not \perp q_{M'}$.
If $q \not \perp r$, then $q_{M'} \not \perp r_{M'}$.
By the Parallelism lemma, 
$p \perp r$ if and only if $p \perp r_{M'}$
and $q \perp r$ if and only if $q_{M'} \perp r_{M'}$.
The conclusion follows from the equivalence axiom.
\end{proof}

\begin{lemma} [Primary base]\label{c4:a:regfin}
If $M'$ is a primary model over $\bigcup_{\eta \in I} M_\eta$,
where $\langle M_\eta \mid \eta \in I \rangle$ is an independent
system and let $p \in S(M')$ be regular.
Then there exists a finite subtree $J \subseteq I$
and a model $M^*$ primary over $\bigcup_{\eta \in J} M_\eta$
such that $p$ is based on $M^*$.
\end{lemma}
\begin{proof}
We prove this by induction on $|I|$.
Let $(J_i : i < |I|)$ be induced subtrees of $I$ of size less than $|I|$
such that $\bigcup_{i < |I|} J_i = I$.
Choose $M_i$ primary over $\bigcup_{\eta \in J_i} M_{\eta}$ such
that $(M_i : i < |I|)$ is $\prec$-increasing, continuous and 
$M' = \bigcup_{ i < |I|} M_i$.
This is possible by existence of primary models over independent system
and the Uniqueness of primary models axiom.
By Local Character, there exists $i < |I|$ such that $p$ is based on $M_i$.
Since the size of $J_i$ is less than $|I|$, we are done by induction. 
\end{proof}

\begin{lemma} \label{c4:l:primeperp}
Let $p=\gatp(a / M)$ be regular and suppose that $p \perp M_1$,
with $M_1 \subseteq M$.
Then $M(a)/M \perp M_1$.
\end{lemma}
\begin{proof}
By axiom (Perp I) it is enough to show that any regular type
$q \in \gaS(M)$ realized in $M(a) \setminus  M$ is orthogonal to any regular 
type $r$ over $M_1$.
But, if $q$ is regular realized in $M(a) \setminus  M$, 
then by
Equivalence we must have $q \not \perp p$.
Since $p \perp M_1$ by assumption, then $p \perp r$.
Then, by definition, $q \perp r$ if and only if
$q \perp r_M$.
Hence, we conclude by Equivalence.
\end{proof}

\begin{lemma}
Let $M_0 \subseteq M$ and let $\bar{a}_1, \bar{a}_2$ such
that $\bar{a}_1 \nonfork_{M} \bar{a}_2$.
Suppose that $\gatp(\bar{a}_\ell / M) \perp M_0$, for $\ell=1,2$.
Let $B$ be such that $B \nonfork_{M} M_0$,
then $\bar{a}_1 \bar{a}_2 \nonfork_M B$.
\end{lemma}
\begin{proof}
By finite character of independence, it is enough
to prove this for finite $B$.
Let $\bar{b}$ be finite such that
\begin{equation} \tag{*}
\bar{b} \nonfork_{M} M_0,
\end{equation}
First, since $\gatp(\bar{a}_2 / M) \perp M_0$, (*) implies that
\begin{equation} \tag{**}
\bar{a}_2 \nonfork_{M} \bar{b}.
\end{equation} 
Thus, by symmetry, we must have $\bar{b}\nonfork_{M} \bar{a}_2$,
so $\bar{b}\nonfork_{M} M_2$, where $M_2$ is primary over $M \cup \bar{a}_2$.
This shows that $\bar{b} \models \gatp(\bar{b}/M)_{M_2}$.
By assumption, we have that
\begin{equation} \tag{***}
\bar{a}_1 \nonfork_M \bar{a}_2, 
\end{equation}
so $\bar{a}_1 \nonfork_M M_2$.
Then 
$\bar{a}_1 \models \gatp(\bar{a}_1 /M)_{M_2}$.
But, $\gatp(\bar{a}_1/M) \perp \gatp(\bar{b}/M)$, so by definition,
we must have
$\bar{a}_1 \nonfork_{M_2} \bar{b}$.
By the first axiom of the independence relation, 
we have $\bar{a}_1 \nonfork_{M_2} \bar{b}\bar{a}_2$.
By transitivity (and dominance) using (***), we obtain
$\bar{a}_1 \nonfork_M \bar{b}\bar{a}_2$.
Hence, by the concatenation property of independence
and ($\dagger$) again, 
we can derive
\[
\bar{a}_1 \bar{a}_2 \nonfork_{M_0} \bar{b}, 
\]
which is what we wanted.
\end{proof}

\begin{corollary}\label{ortho}
Let $M \subseteq N$.
Let $\langle A_i \mid i < \alpha \rangle$ be independent over $N$,
such that $A_i/N \perp M$, for each $i < \alpha$.
Let $B$ be such that $B \nonfork_M N$.
Then $\bigcup \{ A_i \mid i < \alpha \} \nonfork_M B$.
\end{corollary}
\begin{proof}
By finite character of independence and monotonicity, 
we may assume that $\alpha < \omega$.
We prove the statement by induction on $\alpha$ and
use the previous lemma at the successor step.
\end{proof}

\begin{corollary} \label{c4:l:independent}
Let $\langle M_\eta\mid \eta \in I \rangle$ be a system
satisfying:
\begin{enumerate}
\item 
$\langle M_\eta \mid \eta^- =\nu, \eta \in I \rangle$ 
is independent over $M_\nu$, for every $\nu \in I$;

\item 
The type $\gatp(M_\eta / M_{\eta^-}) \perp M_{\eta^{--}}$,
for every $\eta \in I$.
\end{enumerate}
Then $\langle M_\eta \mid \eta \in I \rangle$ is an independent system.
\end{corollary}  

\begin{proof} 
By the finite character of independence, 
we may assume that $I$ is finite.
We prove this statement by induction on $|I|$.
First, notice that if there is no $\eta \in I$
such that $\eta^{--}$ exists, then the result
follows from (1).
We must show that
\[
M_\eta \nonfork_{M_{\eta^-}} 
\bigcup \{M_\nu \mid \eta \not \prec \nu, \nu \in I\}.
\]
Choose $\nu \in I$ of maximal length such that $\eta \not \prec \nu$.
Let 
\[
I_1:= \{ \rho \mid \eta \not \prec \rho, \nu^{-} \prec \rho
\text{ and } \rho \not = \nu \}.
\]
Then, by (1), the system
\begin{equation} \tag{*}
\langle M_\rho, M_\nu \mid \rho \in I_1 \rangle,
\text{ is independent over $M_{\nu^-}$}.
\end{equation}
Let
\[
I_2:= \{ \rho \mid \eta \not \prec \rho, \nu^{-} \not \prec \rho \}.
\]
By induction hypothesis
\begin{equation} \tag{**}
M_{\nu^-} \nonfork_{M_{\nu^{--}}} M_{I_2}M_{\eta}.
\end{equation}
Hence, by the previous corollary, using (*), symmetry on (**) and the
fact that $M_\rho/M_{\nu^-} \perp M_{\nu^{--}}$,
for $\rho \in I_1$ or $\rho = \nu$, we conclude that
\begin{equation} \tag{***}
M_\nu M_{I_1} \nonfork_{M_{\nu^-}} M_{I_2} M_{\eta}.
\end{equation}
Now, by induction hypothesis, we must have 
$M_{I_2} \nonfork_{M_{\nu^-}} M_{\eta}$,
so by concatenation, we must have
\begin{equation} \tag{$\dagger$}
M_\eta \nonfork_{M_{\nu^-}} M_{I_1} M_{I_2} M_{\nu}.
\end{equation}
Now, $M_\eta \nonfork_{M_{\eta^-}} M_{\nu^-}$ by monotonicity and
induction hypothesis.
Therefore, using ($\dagger$),
transitivity and the definition of $I_1$ and $I_2$, 
we conclude that
\[
M_\eta 
\nonfork_{M_{\eta^-}} \bigcup 
\{ M_\nu \mid \eta \not \prec \nu, \nu \in I \}.
\]
\end{proof}

We now come to the main definition of this section.

\begin{definition} 
$\K$ has NDOP if
for every $M_0, M_1, M_2 \in \K$ such that 
$M_1 \nonfork_{M_0} M_2$, for every $M'$ primary over 
$M_1 \cup M_2$
and for every regular type $p \in \gaS(M')$.
Either $p \not \perp M_1$
or $p \not \perp M_2$. 
\end{definition}

\begin{theorem} \label{c4:t:ndop}
Suppose $\K$ has NDOP. 
Let $M, M_\eta \in \K$, for $\eta \in I$ be such that
$\langle M_\eta \mid \eta \in I\rangle$
is an independent system and $M$ is primary over it.
Let $a \in \mathfrak{C} \setminus M$ be such that
$\gatp(a/M)$ is regular.
Then there is $\eta$ such that $\gatp(a/M) \not \perp M_\eta$.
\end{theorem}

\begin{proof}
Let $p = \gatp(a/M)$.
Suppose that $p \perp M_\eta$ for every $\eta \in I$.
By the prime base axiom and parallelism 
we may assume that $I$ is finite.
We will obtain a contradiction to NDOP by induction
on $|I|$.

If $I = \{ \eta \restriction k : k < n \}$, it is obvious
because $\bigcup_{\nu \in I} M_\nu = M_\eta$, 
so by definition of prime, we have $M' = M_\eta$.
But $p \not \perp p$ 
by triviality of independence.
Therefore, $p \not \perp M_\eta$ by definition.

Otherwise, there exists $\nu \in I$ such that both subtrees 
$I_1 := \{ \eta : \eta \in I \ \nu \prec \eta \}$
and $I_2 := \{ \eta : \eta \in I \ \nu \not \prec \eta \}$
are nonempty.
By the third axiom on prime models,
we can choose
$M_k$ prime over $\bigcup_{\eta \in I_k} M_\eta$ for $k = 1, 2$.
By induction hypothesis, we have
\[
p \perp M_1
\qquad
\text{and}
\qquad
p \perp M_2.
\]
Furthermore, since 
$\langle M_\eta \mid \eta \in I\rangle$ is an independent system, 
we have 
\[
\bigcup_{\eta \in I_1} M_\eta \nonfork_{M_\nu} 
\bigcup_{\eta \in I_2} M_\eta.
\]
Therefore, by the symmetry of independence and dominance, we must have
\[
M_1 \nonfork_{M_\nu} M_2.
\]
But, $M'$ is primary over $M_1 \cup M_2$.
This contradicts the fact that $\K$ has NDOP.
\end{proof}

An $\omega$-tree is simply a tree of height at most $\omega$.

\begin{definition} \label{c4:d:decomp}
We say that $\langle M_\eta, a_\eta \mid \eta \in I\rangle$ 
is a \emph{decomposition of $N$ over $M$} if
it satisfies the following conditions:

\begin{enumerate}
\item 
$I$ is an $\omega$-tree;

\item
$\langle M_\eta \mid \eta \in I \rangle$ is a system with 
$M_\eta \subseteq N$ for each $\eta \in I$;

\item 
If $\eta^{--}$ exists for $\eta \in I$, then
$M_\eta / M_{\eta^-} \perp M_{\eta^{--}}$;

\item
For every $\nu \in I$
the system $\langle M_\eta \mid \eta^- = \nu, \eta \in I \rangle$
is independent over $M_\nu$.

\item
$M_{\langle \rangle} = M$ and
$M_\eta$ is primary over $M_{\eta^-} \cup a_\eta$;

\item
For every $\eta \in I$, 
the type $\gatp(a_\eta / M_{\eta^-})$ is regular.
\end{enumerate}
We say that $\langle M_\eta, a_\eta \mid \eta \in I\rangle$ 
is a \emph{decomposition of $N$} if it is a decomposition
of $N$ over $M_{\langle \rangle}$ the primary model over
the empty set.
\end{definition}

Fix $N \in \K$ and $M \prec N$.
We can introduce an ordering between decompositions of $N$ over $M$ 
as follows:
We say that 
\[
\langle M_\eta, a_\eta \mid \eta \in I\rangle
\prec   
\langle N_\eta, b_\eta \mid \eta \in J\rangle
\]
if $I \subseteq J$ and for every $\eta \in I$ we have
\[
M_\eta = N_\eta,
\quad
\text{and}
\quad
a_\eta = b_\eta.
\] 

It is now easy to show that
the set of decompositions of $N$ is inductive:
Let $\langle \mathcal{S}_i \mid i < \alpha \rangle$ be a chain
of decompositions 
$\mathcal{S}_i = \langle M^i_\eta, a^i_\eta \mid \eta \in I^i \rangle$.
First, let $I := \bigcup_{i<\alpha} I^i$.
Then $I$ is an $\omega$-tree.
Hence, we can define the system
$\mathcal{S}:= \langle M_\eta, a_\eta \mid \eta \in I \rangle$,
by $M_\eta := M^i_\eta$ if $\eta \in I^i$ and
$a_\eta := a^i_\eta$, if $\eta \in I^i$.
This is well-defined since 
$\langle \mathcal{S}_i \mid i < \alpha \rangle$
is chain.
We need to check that $\mathcal{S}$ is a decomposition of $N$.
The only nontrivial fact is to check that
for every $\nu \in I$
the system 
\[
\langle M_\eta \mid \eta^- = \nu, \eta \in I \rangle
\]
is independent over $M_\nu$.
If it failed, then by finite character,
there would be a finite set $F \subseteq I$ such that
\[
\langle M_\eta \mid \eta^- = \nu, \eta \in F \rangle
\]
is not independent.
By then, there exists $i < \alpha$ such
that $F \subseteq I^i$, contradicting
the fact that $\mathcal{S}_i$ is
a decomposition of $N$.

Recall that we say that a model $N$ is \emph{minimal} over $A$
if primary models exist over $A$ and if $M(A) \subseteq N$
is primary over $A$, then $N= M(A)$.
Note that a decomposition as in the next theorem is called
\emph{complete}.

\begin{theorem} \label{dec}
Suppose $\K$ has NDOP.
Then for every $M \subseteq N$,
there exists $\langle M_\eta, a_\eta \mid \eta \in I \rangle$
a decomposition of $N$ over $M$ 
such that 
$N$ is primary and minimal over $\bigcup_{\eta \in I} M_\eta$.
\end{theorem}

\begin{proof}
First, notice that the set of decompositions of $N$ over $M$
is not empty.
Therefore, by Zorn's Lemma, since
the set of decompositions of $N$ over $M$ is inductive,
there exists a maximal decomposition 
\begin{equation} \tag{*}
\langle M_\eta, a_\eta \mid \eta \in I \rangle.
\end{equation}
By Lemma \ref{c4:l:independent}, we know that
$\langle M_\eta \mid \eta \in I \rangle$ is
an independent system.
Therefore, by the third axiom for primary models,
there exists $M' \subseteq N$ primary over
$\Union_{\eta \in I} M_\eta$.
We will show that $M'=N$. 
This will show that $N$ is 
primary and minimal over 
$\Union_{\eta \in I} M_\eta$.

Suppose that $M' \not = N$. 
Then, by the axiom of existence of regular types,
there exists a regular type $p \in \gaS(M')$ realized in $N \setminus  M'$.
We are going to contradict the maximality of 
$\langle M_\eta, a_\eta \mid \eta \in I \rangle$.
Since $\K$ has NDOP, by Theorem~\ref{c4:t:ndop}, 
there exists $\eta \in I$ such that $p \not \perp M_\eta$.
Choose $\eta$ of smallest length such that $p \not \perp M_\eta$.
By axiom (Perp I), there exists a regular type $q \in \gaS(M_\eta)$
such that $p \not \perp q$.
Since $q$ is stationary, 
we can choose $q_{M'}$ the unique free extension
of $q$ to the prime model $M'$.
Then, by Lemma~\ref{c4:l:regular}, the type $q_{M'}$
is regular.
Since $p \not \perp q$ and $p \in \gaS(M')$, 
by definition $p \not \perp q_{M'}$.
By Equivalence, since $p$ is realized in $N \setminus M'$, 
there exists $a \in N\setminus M'$ realizing
$q \restriction M'$. 
Hence $\gatp(a/M') = q_{M'}$ and by choice of $q_{M'}$, 
this implies that 
\begin{equation} \tag{**}
a \nonfork_{M_\eta} M'.
\end{equation}
Since $\gatp(a / M_\eta)$ is regular and $a \in N \setminus  M_\eta$, by
the second axiom on primary models,
there exists a primary model $M(a) \subseteq N$
over $M_\eta \cup a$.
By dominance and (**) we must have
\[ 
M(a) \nonfork_{M_\eta} M'.
\]
Thus, by monotonicity of independence and choice of $M'$, 
we conclude that
\begin{equation} \tag{***}
M(a) \nonfork_{M_\eta} \bigcup \{ M_\nu \mid \nu^- = \eta,\  \nu \in I \}.
\end{equation}
But $\{ M_\nu \mid \nu^- = \eta \}$ is independent by definition
of decomposition.
Thus, (***) and Lemma~\ref{c4:l:indepfamily} implies
that 
\[
\{ M_\nu , M(a) \mid \nu^- = \eta, \ \nu \in I \}
\]
is independent over $M_\eta$.
Suppose now that $\eta^-$ exists.
By choice of $\eta$ we must have $p \perp M_{\eta^-}$.
Since $p \not \perp \gatp(a / M_\eta)$, we must
have by Lemma~\ref{c4:l:equiv} and axiom (Perp I) that
$\gatp(a / M_\eta) \perp M_{\eta^-}$.
Hence, by Lemma \ref{c4:l:primeperp}, 
we must have $M(a) / M_\eta \perp M_{\eta^-}$.
This shows that we can add $a/M_{\eta})$ and $M(a)$ to (*) 
and still have
a decomposition of $N$.
This contradicts the maximality of (*).
Thus $N$ is primary and minimal over $\bigcup_{\eta \in I} M_\eta$.
\end{proof}

\begin{corollary}
If $\K$ has NDOP and $N \in \K$, then
there exists a complete decomposition.
\end{corollary}
\begin{proof}
By the previous theorem since by axiom on primary models there
exists a primary model over the empty set.
\end{proof}

The same proof shows:

\begin{corollary}
If $\K$ has NDOP and $N \in \K$ is primary over a decomposition
$\langle M_\eta \mid \eta \in I \rangle$ of $N$ over $M$,
then $\langle M_\eta \mid \eta \in I \rangle$ is a complete
decomposition of $N$ over $M$.
\end{corollary}

\subsection{Examples}

The abstract decomposition given in this section generalises
the known NDOP cases.

There are several classical first order cases.
The first one is for $\aleph_0$-saturated models of a totally transcendental
theory $T$. 
A second one is for
$\aleph_\epsilon$-saturated models of a superstable theory $T$.
And finally, for the class of models of a totally transcendental theory $T$.
In each case, $\mathfrak{C}$ can be taken to be the saturated
monster model for $T$.
The {\em independence relation} is nonforking. 
{\em Regular types} in the first two cases are just the regular types in the
sense of first order. 
In the last case, they correspond to strongly regular types.
The {\em primary models} are the $\mathbb{F}^s_{\aleph_0}$-primary models,
the $\aleph_\epsilon$-primary models (also called $\mathbb{F}^a_{\aleph_0}$)
for the second case, and
the $\mathbb{F}^t_{\aleph_0}$-primary models in the third case.
All the results needed to apply the theorem can be found in
\cite{Sh b}.

In the nonelementary case, there are two published examples.
One in the context of 
an excellent Scott sentence in $L_{\omega_1, \omega}$
\cite{GrHa}. 
Excellence implies AP and JEP;
The model $\mathfrak{C}$ can be taken to be any sufficiently large
full model.
The {\em independence relation} is that afforded by the rank.
Regular types are the SR types.
The {\em primary models} are the usual primary models
and their uniqueness is clear.
The existence of primary models follows from excellence
(see \cite{Sh:871}, \cite{Sh:872}) and
the relevant orthogonality calculus can be found in~\cite{GrHa}).
The other is for the class of locally saturated models of 
a superstable diagram.
There we have amalgamation over sets, so $\mathfrak{C}$ can be 
taken to be strongly homogeneous. 
The details are in \cite{HySh:2}.

The aim of the next section is to prove that the axiomatic framework
developed in this section holds 
for the class $\mathcal{K}$ of 
$(D,\aleph_0)$-homogeneous models of a totally transcendental $D$.
As we pointed out in the previous paragraph, we have an abstract
elementary class with AP and JEP, and even more: we can work inside
a large homogeneous model.

The \emph{independence relation} is given by the rank;
the axioms for independence, 
the existence of stationary types, the existence of regular types
can all be found in \cite{Le1}.
The {\em primary models} are the $D^s_{\aleph_0}$-primary models;
their uniqueness is clear and 
their existence over all sets in the totally transcendental case
is also proved in \cite{Le1}.

This leaves us with the proof of Equivalence, and
Dominance. 
These results are part of what is called Orthogonality Calculus.

\bigskip
\section{Orthogonality calculus in finite diagrams}

In this section, we work in the context of totally transcendental
good diagrams. 

Let $T$ be a complete first order theory in a language $L$.
A type $\ftp(c/A, M)$ is simply the set of first order formulas over $A$
which are true of $c$ in $M$.
A {\em diagram} $D$ is a set of the form $\{ \ftp(c/\emptyset,M) : c \in M \}$
for some $M \models T$.
Fix a diagram $D$.
A set $A$ is a {\em $D$-set}, if $\ftp(c/\emptyset) \in D$ for each
finite $c \in A$.
A $D$-model is a model whose universe is a $D$-set.
We are interested in the class of $D$-models, {\em i.e.} 
the nonelementary class of models of $T$
omitting, over the empty set, all the types outside $D$.

This leads to the following notion of types:
\begin{definition}
$S_D(A) = \{ p \in S(A) : \text{ $Ac$ is a $D$-set for each $c \models p$}\}$.
\end{definition} 

A model $M$ is {\em $\lambda$-homogeneous} if whenever $f : M \rightarrow M$
is a partial elementary map with $| \dom(f)| < \lambda$ and
$a \in M$ then there is an elementary map $g : M \rightarrow M$
extending $f$ such that $a \in \dom(g)$.
A model $M$ is $(D,\lambda)$-homogeneous if $M$ 
is $\lambda$-homogeneous and $M$ realizes exactly the types in $D$.
Then, if $p \in S_D(A)$ and $A \subseteq M$ has size less
than $\lambda$, then $p$ is realized in $M$ if $M$ is 
$(D,\lambda)$-homogeneous.  

We will work inside a large $(D,\bar{\kappa})$-homogenous model $\mathfrak{C}$
of size $\bar{\kappa}$,
which functions as our monster model.
Satisfaction is defined with respect to $\mathfrak{C}$,
and all sets and models are assumed to be inside $\mathfrak{C}$,
so all the relevant types are realized in $\mathfrak{C}$.
The existence of such a model is the meaning of {\em good}.

Here is the meaning of stability in this context:
\begin{definition}
$D$ is {\em $\lambda$-stable} if
$|S_D(A)| \leq \lambda$ for each $D$-set of size~$\lambda$. 
\end{definition}
In \cite{Le1}, a notion of rank is introduced which is shown to
be bounded under $\omega$-stability. 
$D$ is said to be {\em totally transcendental}
if the rank is bounded.
In the rest of this paper, we assume that $D$ is a totally
transcendental diagram.
We already established in \cite{Le1} that many of the axioms of the previous
section hold for totally transcendental diagrams (notably the properties
of the independence relation and the existence of primary models)
and facts from \cite{Le1} will be used freely.
We will now develop what is referred to as \emph{orthogonality calculus}
for this context and show that the remaining axioms used
to obtain an abstract decomposition theorem also hold
for the class of $(D, \aleph_0)$-homogeneous models of a totally
transcendental diagram $D$.

Notice that homogeneity implies that the 
notion of $\gatp(\bar a/M,N)$ coincides with $\ftp(\bar a/M,N)$.

The next few lemmas show Dominance.

First, for $D$-sets $A$ and $B$, we say that 
$A \subseteq_{TV} B$, if every $D$-type over finitely many
parameters in $A$ realized in $B$ is realized in $A$.
The subscript TV stands for Tarski-Vaught.

\begin{lemma} \label{TV}
Let $M$ be $(D, \aleph_0)$-homogeneous.
Suppose $\bar{a} \nonfork_M \bar{b}$.
Then, for every $\bar{m} \in M$
the type $\ftp(\bar{b}/ \bar{m} \bar{a})$ is realized in $M$.
\end{lemma}
\begin{proof}
By symmetry, $\bar{b} \nonfork_M \bar{a}$.
Hence, by taking a larger $\bar{m}$ if necessary,
we may assume that $\ftp(\bar{b}/M\bar{a})$ does not split over $\bar{m}$.
By $(D,\aleph_0)$-homogeneity of $M$, we
can find $\bar{b}' \in M$, such that
$\ftp(\bar{b}/\bar{m}) = \ftp(\bar{b}'/\bar{m})$.
We claim that $\ftp(\bar{b}/\bar{m} \bar{a})=\ftp(\bar{b}'/\bar{m} \bar{a})$.
If not, there exists a formula $\phi(\bar{x}, \bar{m}, \bar{a})$
such that $\models \phi[\bar{b},\bar{m}, \bar{a}]$ and
$\models \neg \phi[\bar{b}',\bar{m}, \bar{a}]$.
But, $\ftp(\bar{b}/\bar{m}) = \ftp(\bar{b}'/\bar{m})$,
so $\ftp(\bar{a}/M\bar{b})$ splits over $\bar{m}$,
a contradiction.
\end{proof}

The next lemma is standard. 
Recall that $p \in S(A)$ is {\em $D^s_{\aleph_0}$-isolated} if
$p \in S_D(A)$ and is $\mathbb{F}^s_{\aleph_0}$-isolated.

\begin{lemma}
Let $A, B$ be $D$-sets such that $A \subseteq_{TV}B$.
If $\ftp(\bar{c}/A)$ is $D^s_{\aleph_0}$-isolated,
then $\ftp(\bar{c}/A) \vdash \ftp(\bar{c}/B)$.
\end{lemma}
\begin{proof}
Let $q(\bar{x},\bar{a}) \vdash \ftp(\bar{c}/A)$, with $\bar{a} \in A$.
Suppose that $\ftp(\bar{c}/A) \not \vdash \ftp(\bar{c}/B)$.
Then, there exists $\bar{b} \in B$ and a formula 
$\phi(\bar{x}, \bar{y})$
such that $q(\bar{x},\bar{a}) \cup \phi(\bar{x}, \bar{b})$
and $q(\bar{x},\bar{a}) \cup \neg \phi(\bar{x}, \bar{b})$
are both realized in $\mathfrak{C}$.
By assumption,
there exists $\bar{b}' \in A$ 
realizing be such that $\ftp(\bar{b}/\bar{a}) = \ftp(\bar{c}/\bar{c})$.
Hence, by an automorphism fixing $\bar{a}$ and sending $\bar{b}$
to $\bar{b}'$, both $q(\bar{x},\bar{a}) \cup \phi(\bar{x}, \bar{b}')$
and $q(\bar{x},\bar{a}) \cup \neg \phi(\bar{x}, \bar{b}')$
are realized in $\mathfrak{C}$.
This contradicts the choice of $q(\bar{x},\bar{a})$.
\end{proof}

A model is $M$ is {\em $D^s_{\aleph_0}$-primary over $A$}
if $M = A \cup \{ a_i : i < \lambda \}$ and 
$\ftp(a_i/A \cup \{ a_j : j < i \})$
is $D^s_{\aleph_0}$-isolated for each $i < \lambda$.
We denote by $M(A)$ the $D^s_{\aleph_0}$-primary
model over $M \cup A$.

\begin{theorem}[Dominance]
Let $M$ be $(D,\aleph_0)$-homogeneous and $A$ be a $D$-set.
For each $B$, if $A \nonfork_M B$, then $M(A) \nonfork_M B$.
\end{theorem}
\begin{proof}
By finite character of independence, it is enough to
show that if $\bar{a} \nonfork_M \bar{b}$, then
$\bar{c} \nonfork_M \bar{b}$, for each finite $\bar{c} \in M(\bar{a})$.
Let $\bar{c} \in M(\bar{a})$ be given.
Then $\ftp(\bar{c}/M\bar{a})$ is $D^s_{\aleph_0}$-isolated.
Hence, by assumption and Lemma \ref{TV}, 
$\ftp(\bar{c}/M\bar{a}) \vdash \ftp(\bar{c}/M\bar{a}\bar{b})$.
Therefore, $\bar{c} \nonfork_M \bar{b}$.
\end{proof}

Recall the definition of orthogonality.

\begin{definition} 
Let $p \in S_D(B)$ and $q \in S_D(A)$ be
stationary. We say that $p$ is \emph{orthogonal} to $q$, written $p \perp q$, if
for every $D$-model $M$ containing $A \cup B$
and for every $a \models p_M$ and $b \models q_M$,
we have $a \nonfork_M b$;
\end{definition}

Recall the following technical lemma in \cite{Le1}.
\begin{lemma} \label{tech} Let 
$p, q \in S_D(M)$ and $M \subseteq N$ be in $\mathcal{K}$.
If
$a \nonfork_M b$
for every $a \models q$ and $b \models p$,
then 
$a \nonfork_N b$ for every $a \models q_N$ and $b \models p_N$.
\end{lemma}

Then, by the previous lemma we can immediately simplify
the definition: 
for $p, q \in S_D(M)$,
we have $p \perp q$ if and only if $\bar{a} \nonfork_M \bar{b}$
for every $\bar{a} \models p$ and $\bar{b} \models q$.

The following lemma is also in \cite{Le1}.

\begin{lemma} \label{halfequiv} Let $M$ be $(D,\aleph_0)$-homogeneous.
If $\bar{a} \nonfork_M \bar{b}$ and $\ftp(\bar{a}/ M \bar{b})$ 
is $D^s_{\aleph_0}$-isolated, then $\bar{a} \in M$.
\end{lemma}

\begin{lemma} \label{otherhalfequiv}
Let $\ftp(\bar{a}/M\bar{b})$ be isolated,
and $\ftp(\bar{b}/M)$ be regular.
Suppose that $\bar{a} \fork_M \bar{b}$.
Then, for any $\bar{c}$ if 
$\bar{a} \nonfork_M \bar{c}$, then
$\bar{b} \nonfork_M \bar{c}$.
\end{lemma}

\begin{proof}
Suppose that $\bar{b} \fork_M \bar{c}$.
By symmetry, we have that $\bar{c} \fork_M \bar{b}$.
Choose $q(\bar{x}, \bar{m}, \bar{b}) \subseteq \ftp(\bar{c}/M\bar{b})$
such that 
\[
R[q(\bar{z}, \bar{m}, \bar{b})]=
R[\ftp(\bar{c}/M\bar{b})] < 
R[\ftp(\bar{c}/M)].
\]
Without loss of generality,
since $\bar{a} \fork_M \bar{b}$, 
we can choose $p(\bar{x}, \bar{m}, \bar{b}) \subseteq \ftp(\bar{a}/M\bar{b})$
be such that 
\[
R[p(\bar{y}, \bar{m}, \bar{a})]=
R[\ftp(\bar{b}/M\bar{a})] < 
R[\ftp(\bar{b}/M)]
\]
and also
\[
R[p(\bar{b}, \bar{m}, \bar{x})]=
R[\ftp(\bar{a}/M\bar{b})] < 
R[\ftp(\bar{a}/M)].
\]
Choose $\bar{c}' \in M$ such that 
$tp(\bar{c}/\bar{m})= \ftp(\bar{c}'/\bar{m})$.
Since $\bar{a} \nonfork_M \bar{c}$, we have in particular
that $\ftp(\bar{a}/M\bar{c})$ does not split over $\bar{m}$
so that $\ftp(\bar{c}/\bar{m}\bar{a})=\ftp(\bar{c}'/\bar{m}\bar{a})$.
Thus, $\bar{b}$ realizes the following type
\begin{equation} \tag{*}
p(\bar{y}, \bar{m}, \bar{a}) \cup 
q(\bar{a}, \bar{m}, \bar{y}) \cup
\ftp(\bar{b}/\bar{m}).
\end{equation}
Since $\ftp(\bar{a}/M\bar{b})$ is isolated, we may assume that
$M(\bar{a}) \subseteq M(\bar{b})$.
Now choose $\bar{b}' \in M(\bar{a})$ realizing (*).
If $\bar{b}' \in M$, then
$R[\ftp(\bar{a}/M)] \leq R[\ftp(\bar{a}/\bar{m}\bar{b}')]
= R[p(\bar{b}',\bar{m},\bar{x})]$, a contradiction.
Hence $\bar{b}' \not \in M$ and so
$\bar{b}' \fork_M \bar{b}$, by the previous lemma.
Thus $\ftp(\bar{b}'/M)$ extends $\ftp(\bar{b}/\bar{m})$
and is not orthogonal to it, thus
since $\ftp(\bar{b}/M)$ is regular based on $\bar{m}$,
we must have $\ftp(\bar{b}'/M)= \ftp(\bar{b}/M)$.
This is a contradiction, since then,
$\bar{b}'$ realizes $q(\bar{c}', \bar{m}, \bar{y})$.
\end{proof}

The next corollary is Equivalence.

\begin{corollary} [Equivalence]
Let $M \in \mathcal{K}$,
let $p, q \in S_D(M)$ be regular,
and let $\bar{b} \not \in M$ realize $p$. 
Then $q$ is realized in $M(\bar{b}) \setminus M$ if and only
if $p \not \perp q$.
\end{corollary}

\begin{proof} 
Let $\bar{b} \in M$ realize $p$.
Let $M(\bar{b})$ be $D^s_{\aleph_0}$-primary over $M \cup \bar{b}$.

Let $\bar{a} \in M(\bar{b}) \setminus M$.
Then $\ftp(\bar{a}/M\bar{b})$ is $D^s_{\aleph_0}$-isolated.
If $p \perp q$, then $\bar{b} \nonfork_M \bar{a}$.
Hence, by symmetry $\bar{a} \nonfork_M \bar{b}$, and
so $\bar{a} \in M$, by Lemma \ref{halfequiv},
a contradiction.

For the converse, suppose that $p \not \perp q$.
This implies that there is $\bar{a} \models q$ such that
\[
\bar{a} \fork_M \bar{b}.
\]
Let $q(\bar{x}, \bar{m}, \bar{b}) \subseteq \ftp(\bar{a}/M\bar{b})$
be such that 
\[
R[q(\bar{x}, \bar{m}, \bar{b})] = R[\ftp(\bar{a}/M\bar{b})] < R[ q].
\]
Since $q$ is regular, we may further assume that $q$ is based on $\bar{m}$.
Thus, the element $\bar{a}$ realizes the type
\begin{equation} \tag{*}
q(\bar{x}, \bar{m}, \bar{b}) \cup q \restriction \bar{m}.
\end{equation}
Since $M(\bar{b})$ is in particular $(D, \aleph_0)$-homogeneous,
there is $\bar{a}' \in M(\bar{b})$ realizing the type (*).
Since $M(\bar{b})$ is $(D, \aleph_0)$-primary, we must have
that $\ftp(\bar{a}'/M\bar{b})$ is isolated.
Thus, since $\bar{b} \fork_M \bar{a}$, we must
have by the Lemma \ref{otherhalfequiv} that $\bar{a}' \fork_M \bar{a}$.
This implies that $\ftp(\bar{a}' /M )$ is an extension of 
the regular type $q \restriction \bar{m}$ which is not
orthogonal to $q$. 
Hence, since $q$ is regular, we must have $q = \ftp(\bar{a}' / M)$.
This shows that $q$ is realized (by $\bar{a}'$) in $M(\bar{b})$.
\end{proof}

We encountered Morley sequences when we talked about
stationary types in the previous chapter.
The definition can be made for any type.

\begin{definition}
Let $p \in S_D(A)$.
We say that $\langle \bar{a}_i \mid i < \omega \rangle$
is a \emph{Morley sequence} for $p$ if
\begin{enumerate}
\item
The sequence
$\langle \bar{a}_i \mid i < \omega \rangle$ is indiscernible over $A$;

\item
For every $i < \omega$ we have 
$\bar{a}_i \nonfork_A A \cup \{\bar{a}_j \mid j < i \}$.
\end{enumerate}
\end{definition}

The next fact was established in the previous chapter.

\begin{fact} 
If $p \in S_D(A)$ is stationary,
then there is a Morley sequence for~$p$.
\end{fact}

The next theorem is Axiom (Perp I).
\begin{theorem}[Perp I]
Let $p \in S_D(N)$ be regular, $M \subseteq N$.
Then $p \perp M$ if and only if $p \perp q$, for every regular 
$q \in S(M)$.
\end{theorem}
\begin{proof}
One direction is obvious.
Suppose that $p \not \perp M$.
We will find a regular type $q \in S_D(M)$ such that $p \not \perp q$.

Since $p$ is regular, there
exists a finite set $\bar{f} \subseteq N$ such that
$p$ is regular over $\bar{f}$. 
Write $p(\bar{x}, \bar{f})$ for the stationary type $p_{\bar{f}}$.
Also, there exists a finite set $\bar{e} \subseteq M$
such that $\ftp(\bar{f}/M)$ is based on $\bar{e}$.
Since $p \not \perp M$, there exists a stationary type
$r \in S(M)$ such that $p \not \perp r$. 
By monotonicity, we can find
$\bar{a} \models p$, $\bar{b} \models r_N$ such
that $\bar{a} \fork_{\bar{e}} \bar{f}\bar{b}$.

Since $M$ is $(D, \aleph_0)$-homogeneous, there
exists $\langle \bar{f}_i \mid i < \omega \rangle \subseteq M$,
a Morley sequence for $\ftp(\bar{f}/\bar{e})$.
Let $p_i : = p(\bar{x},\bar{f}_i)_M$. 
This is well-defined since $p(\bar{x},\bar{f})$ is stationary
and $\ftp(\bar{f}/\bar{e})= \ftp(\bar{f}_i/\bar{e})$, so
$p(\bar{x},\bar{f}_i)$ is stationary.

For each $i < \omega$, we can choose $M_i \subseteq N$
such that there is an automorphism $g_i$ with
$g_i(\bar{f})=\bar{f}_i$, $g_i(\bar{c})=\bar{c}$ and 
$g_i(M)=M_i$. 
Since
$p_{M_i}$ is regular and $p_i = g^{-1}(p_{M_i})$,
then 
\begin{equation} \tag{*}
p_i
\text{ is regular, for each $i < \omega$.}
\end{equation}
A similar reasoning using an automorphisms sending
$\bar{f}\bar{f}_0$ to $\bar{f}_i\bar{f}_j$ shows
that 
\begin{equation} \tag{**}
p \perp p_0
\text{ implies }
p_i \perp p_j,
\text{ for every $i \not = j < \omega$}.
\end{equation}
Finally, using the fact that $p \not \perp r$,
we can derive
\begin{equation} \tag{***}
p_i \not \perp r,
\text{ for every $i<\omega$}.
\end{equation} 
If we show that $p \not \perp p_0$, then
(*) implies the conclusion of the lemma.
Suppose, for a contradiction, that $p \perp p_0$.
By (***) 
we can find $\bar{b}' \models r$ and $\bar{a}_i \models p_i$, such
that $\bar{b}' \fork_M \bar{a}_i$ and $\bar{a}_i \not \in M$,
for each $i < \omega$.
Now (**) implies that 
$\bar{a}_{j+1} \nonfork_M  \{ \bar{a}_i \mid i \leq j \}$,
for every $j < \omega$.
Hence, by (*) and Lemma \ref{halfequiv}, we have
$\bar{a}_{i+1} \not \in M_i$, where $M_i$ is 
$D^s_{\aleph_0}$-primary over $M \cup \{ \bar{a}_j \mid j < i \}$.
Let $N$ be $D^s_{\aleph_0}$-primary over 
$M \cup \{ \bar{a}_j \mid j < \omega \}$.
Since $\kappa(D) = \aleph_0$, there exists $n < \omega$ such
that $\bar{b}' \nonfork_{M_n} N$.
Hence, by monotonicity, $\bar{b}' \nonfork_{M_n} \bar{a}_n$.
By symmetry over models, $\bar{a}_n \nonfork_{M_n} \bar{b}'$.
But $\bar{a}_n \nonfork_M \{ \bar{a}_i \mid i < n\}$,
and so $\bar{a}_n \nonfork_M M_n$, by dominance and symmetry.
Hence, by transitivity of the independence relation, 
we have $\bar{a}_n \nonfork_M \bar{b}'$, so
$\bar{b}' \nonfork_M \bar{a}_n$, a contradiction.
\end{proof}

We now prove two additional lemmas that will be used
in the next section.
\begin{lemma}\label{reguperp}
If $p \in S(M_1)$ is regular,
$p \perp M_0$, and $M_1 \nonfork_{M_0} M_2$, then $p \perp M_2$.
\end{lemma}

\begin{proof}
Suppose that $p \not \perp M_2$.
Then, by definition, 
there exists $q \in S(M_2)$ such that $p \not \perp q$.
By definition, there is $N \supseteq M_1 \cup M_2$
such that 
\begin{equation} \tag{*}
p_N \not \perp q_N.
\end{equation}
We are going to find a type $q' \in S(M_0)$ such that $p \not \perp q'$.

Since $p$ and $q$ are stationary,
there exist finite sets $\bar{c} \subseteq M_1$, 
$\bar{d} \subseteq M_2$, and $\bar{e} \subseteq M_0$ such that
$p$ is based on $\bar{c}$, 
$q$ is based on $\bar{d}$, and
both $\ftp(\bar{c}/M_0)$ and
$\ftp(\bar{d}/M_0)$ are based on $\bar{e}$.

By (*) and finite character, 
there exist a set $F \subseteq N$,
and $\bar{a}$, $\bar{b}$ such
that 
\begin{equation} \tag{**}
\bar{a} \models p_{M_1 M_2 F},
\quad
\bar{b}  \models q_{M_1 M_2 F},
\quad
\text{but}
\quad
\bar{a} \fork_{M_1 M_2 F} \bar{b}.
\end{equation}
By monotonicity, we may assume that
$\bar{c}\bar{d}\bar{e} \subseteq F$.
Since $\ftp(\bar{a}\bar{b}/N)$ is stationary,
we may also assume that $\ftp(\bar{a}\bar{b}/M_1M_2F)$ is stationary based
on $F$. 
Finally, we may further assume that 
$R[\ftp(\bar{a}/\bar{c})] < R[\ftp(\bar{a}/\bar{c}F)]$.

Since $M_0$ is $(D, \aleph_0)$-homogeneous, 
we can choose $\bar{d}' \in M_0$ such that
$\ftp(\bar{d}'/\bar{e})=\ftp(\bar{d}/\bar{e})$.
By stationarity, 
we have $\ftp(\bar{c}\bar{d}\bar{e}/\emptyset)=
\ftp(\bar{c}\bar{d}\bar{e}/\emptyset)$.
Now choose $F' \subseteq M_1$ such that 
$\ftp(\bar{c}\bar{d}\bar{e}F /\emptyset) = 
\ftp(\bar{c}\bar{d}'\bar{e}F /\emptyset)$.
Finally, let $\bar{a}'\bar{b}' \in \mathfrak{C}$
such that 
$\ftp(\bar{a}\bar{b}\bar{c}\bar{d}\bar{e}F /\emptyset)=
\ftp(\bar{a}'\bar{b}'\bar{c}\bar{d}'\bar{e}F' /\emptyset)$.

By invariance under automorphism, we have
$R[\ftp(\bar{a}'/\bar{c})]=R[\ftp(\bar{a}'/F')]$ and
$R[\ftp(\bar{b}'/\bar{d}')]= R[\ftp(\bar{b}'/F')]$,
since these statements are true without the apostrophe $'$.

Now let $q':= \ftp(\bar{b}'/\bar{d}')_{M_0} \in S(M_0)$.
Such a type exists since $\ftp(\bar{b}'/\bar{d}')$ is stationary.
We claim that $p \not \perp q'$.
Otherwise, by the previous remark, we have $p \perp q'_{M_1}$.
Now, let $\bar{a}''\bar{b}'' \models \ftp(\bar{a}'\bar{b}'/F')$.
We have $\bar{a}'' \models p$, $\bar{b}'' \models q'_{M_1}$
and so $\bar{a}'' \models p_{M_1 \bar{b}''}$.
But then 
$R[\ftp(\bar{a}''/\bar{c})]= R[\ftp(\bar{a}''/\bar{b}''F')]$
This contradicts the fact that 
$\ftp(\bar{a}'\bar{b}'/F')= \ftp(\bar{a}''\bar{b}''/F')$.
\end{proof}

\begin{lemma} \label{strongequiv}
Let $p, q \in S_D(M)$ be regular.
Let $\bar{a} \not \in M$ realize $p$.
If $p \not \perp q$, then there exists $\bar{b} \in M(\bar{a})\setminus M$
realizing $q$ such that $M(\bar{a})=M(\bar{b})$.
\end{lemma}
\begin{proof}
By equivalence, there exists $\bar{b} \in M(\bar{a}) \setminus M$
realizing $q$.
By definition of prime, it is enough to show that
$\ftp(\bar{a}/M\bar{b})$ is $D^s_{\aleph_0}$-isolated.

Let $\bar{c} \in M$ be finite such that $p$ is regular over $\bar{c}$,
and write $p(\bar{x},\bar{c}) = p \restriction \bar{c}$.
Now, since $\ftp(\bar{b}/M\bar{a})$ is $D^s_{\aleph_0}$-isolated,
there exists $r_1(\bar{y}, \bar{a})$ over $M$ isolating 
$\ftp(\bar{b}/M\bar{a})$.
By a previous lemma, we know that $\bar{a} \fork_M \bar{b}$, so
let $r_2(\bar{x}, \bar{b})$ witness this.
We claim that the following type isolates $\ftp(\bar{a}/M\bar{b})$:
\begin{equation} \tag{*}
p(\bar{x},\bar{c}) \cup r_1(\bar{b}, \bar{x}) \cup r_2(\bar{x}, \bar{b}).
\end{equation}
Let $\bar{a}' \in M(\bar{a})$ realize (*).
Then, $\bar{a}' \not \in M$ by choice of $r_1$.
Hence, $\bar{a} \fork_M \bar{a}'$ so by choice of $p(\bar{x},\bar{c})$,
we have $\ftp(\bar{a}'/M) = \ftp(\bar{a}/M)$.
Thus, $\ftp(\bar{a}/M\bar{b}) = \ftp(\bar{a}'/M\bar{b})$ 
using $r_2(\bar{a}', \bar{y})$.
\end{proof}

We can now show using the language of Section 1.

\begin{theorem}
Let $\mathcal{K}$ be 
the class of $(D,\aleph_0)$-homogeneous models of a totally
transcendental diagram $D$.
Let $N \in \mathcal{K}$ have NDOP.
Then $N$ has a complete decomposition.
\end{theorem}
\begin{proof}
All the axioms of Section 1 have been checked for $\mathcal{K}$.
\end{proof}

\begin{remark}
Similarly to the methods developed in this section for the class 
of $(D, \aleph_0)$-homogeneous models of a totally transcendental diagram
$D$, we can check all the axioms for the class of 
$(D,\mu)$-homogeneous models of a totally transcendental
diagram $D$, for any infinite $\mu$.
This implies that if $\mathcal{K}$ is the class of 
$(D,\mu)$-homogeneous models of a totally transcendental
diagram $D$ and if $N \in \mathcal{K}$ has NDOP, then
$N$ has a complete decomposition (in terms of models
of $\mathcal{K}$).
\end{remark}

\bigskip
\section{DOP in finite diagrams}

Let $\mathcal{K}$ be the class of $(D,\aleph_0)$-homogeneous models
of a totally transcendental diagram.
In the language of the axiomatic framework, we take
$\mathfrak{C}$ to be a large homogeneous model.
We say that $\mathcal{K}$ satisfies DOP if $\mathcal{K}$
does not have NDOP.
Recall that $\lambda(D)$, the first stability cardinal, is $|D|+|T|$.

\begin{claim}
Suppose that $\mathcal{K}$ has DOP.
Then there exists $M, M_i, M' \in \mathcal{K}$ for $i=1,2$ 
such that 
\begin{enumerate}
\item
$M_1 \nonfork_M M_2$;

\item
$M'$ is prime over $M_1 \cup M_2$;

\item
$\| M' \| = \lambda(D)$;

\item
$M_i=M(\bar{a}_i)$, for $i=1,2$;

\item
There exists a regular type $p \in S(M')$
such that $p \perp M_i$, for $i=1,2$;

\item
The type $p$ is based on $\bar{b}$ and $\ftp(\bar{b}/M_1 \cup M_2)$
is isolated over $\bar{a}_1 \bar{a}_2$.
\end{enumerate}

\end{claim}

\begin{proof}
By assumption, there exists $\mathcal{K}$ 
fails to have NDOP.
Then, there exist $M_i\in \mathcal{K}$
for $i \leq 2$ with $M_1 \nonfork_{M_0} M_2$,
there exists $M''$ which is 
$D^s_{\aleph_0}$-primary over $M_1 \cup M_2$ and 
there exists a regular type $p \in S(M'')$ such that
$p \perp M_i$,
for $i=1,2$.

Let $\bar{b} \in M''$ be a finite set such that 
$p$ is based on $\bar{b}$.
Let $\bar{a}_i \in M_i$, for $i=1,2$ be such that
$\ftp(\bar{b}/M_1 \cup M_2)$ is 
$D^s_{\aleph_0}$-isolated over $\bar{a}_1\bar{a}_2$.
Let $M \in \mathcal{K}$, 
$M \subseteq M_0$ of cardinality $\lambda(D)$ be such that
$\bar{a}_1 \nonfork_M M_0$. 
Such a model exists using local character and prime models.
Let $M(\bar{a}_i)$ be prime over $M \cup \bar{a}_i$, for $i=1,2$.
Then, by Dominance, Transitivity, and Monotonicity, we have 
$M(\bar{a}_1) \nonfork_M M(\bar{a}_2)$.
By axiom on prime there exists $M' \subseteq M''$ prime over 
$M(\bar{a}_1) \cup M(\bar{a}_2)$.
We may assume that $B \subseteq M'$.
Let $p' = p \restriction M'$.
Then $p' \in S(M')$ is regular based on $\bar{b}$ and
$p'_{M''} = p$.
It remains to show that $p' \perp M(\bar{a}_i)$, for $i=1,2$.
Let $r \in S(M(\bar{a}_i)$ be regular.
Then $r_{M_i}$ is regular by our axiom.
Furthermore, by definition, $p' \perp r$ if and only if 
$p' \perp r_{M'}$.
By Parallelism, since $M' \subseteq M''$, it is equivalent
to show that $p \perp r_{M''}$.
But, $r_{M_i} \in S(M_i)$ is regular, 
$p \perp M_i$, and $r_{M''} = (r_{M_i})_{M''}$.
Therefore, by choice of $p$ we have $p \perp r_{M''}$,
which finishes the proof.
\end{proof}

Let $\mu > \lambda(D)$ be a cardinal (for
the following construction, we may have $\mu \geq \lambda(D)$,
but the strict inequality is used in the last claim).
Let $\langle M_i \mid i < \mu \rangle$ be independent over a model 
$M \subseteq M_i$.
Suppose that $\|M_i\|=\lambda(D)$.
Let $R \subseteq [\mu]^2$ and suppose that
$M_s = M(M_i \cup M_j)$, for $s=(i,j)$.
Such a model exist for each $s \in [\mu]^2$ by the axioms on prime. 
Then, by Dominance and the axiom on primes, the following system
is independent:
\begin{equation} \tag{*}
\langle M_i \mid i < \mu \rangle
\cup \{M\} \cup \langle M_s \mid s \in R \rangle.
\end{equation}
Hence, there exists a model
$M_R$ prime over $\bigcup_{i < \mu} M_i \cup \bigcup_{s \in R} M_s$.

Let $s=(i,j)$ and
suppose that there exists a regular type $p_s \in S(M_s)$ 
such that $p_s \perp M_i$, $p_s \perp M_j$.
Let $I_s$ be a Morley sequence for $p_s$ of length $\mu$.
(Such a sequence exists since $\mathfrak{C}$ is $(D,\mu^+)$-homogeneous.
Then, by Dominance, definition of a Morley sequence, and axiom on prime,
there exists $N_s= M_s(I_s)$.

The next claim will allows us to choose prime models over complicated
independent systems with some additional properties.

\begin{claim}
The system 
\[
\mathcal{S}_R = \langle M_i \mid i < \mu \rangle
\cup \{M\} \cup \langle N_s \mid s \in R \rangle
\]
is an independent system.
\end{claim}
\begin{proof}
By definition, we must show that
$N_s \nonfork_{M_i} A$, when $A = \bigcup_{t \in R, t \not = s} N_t$.
By finite character, it is enough to show this for $R$ finite.
We prove this by induction on the cardinality of $R$.
When $R$ is empty or has at most one element,
there is nothing to do.
Suppose that $R = \{ s_i \mid i \leq n \} \cup \{ s \}$.
We show that we can replace $M_{s_i}$ by $N_{s_i}$ and
$M_s$ by $N_s$ and still have an independent system.
By (*), it is enough to show that if $M_s \nonfork_{M_{s^-}} A$,
then $N_s \nonfork_{M_{s^-}} A$, for $A = \bigcup_{i \leq n} N_{s_i}$.
Using the axioms of the independence relation,
it is enough to show that $N_s \nonfork_{M_s} A$.
By induction hypothesis, we have
\begin{equation} \tag{**}
N_s \nonfork_{M_s} \bigcup_{i < n} N_{s_i}
\quad
\text{and}
\quad
N_{s_n} \nonfork_{M_{s_n}} \bigcup_{i < n} N_{s_i}.
\end{equation}
Now, either $s \cap s_n$ is empty so $M_s \nonfork_M M_{s_n}$ 
by (*)
or they extend $j$ and so $M_s \nonfork_{M_j} M_{s_n}$,
by (*) again.
Since $M \subseteq M_j$, in either case, $p_s \perp M_j$,
by choice of $p_j$.
Hence $p_s \perp M_{s_n}$ using Lemma \ref{reguperp}.
By induction hypothesis, there exists $N'$ a prime model 
over $\bigcup_{i < n} N_{s_i}$.
Hence, by (**) and Dominance 
$N_s \nonfork_{M_s} N'$ and $N_{s_n} \nonfork_{M_{s_n}} N'$.
Hence, using again by Lemma \ref{reguperp},
we have $p_s \perp N'$.
Thus, $I_s \nonfork_{M_{s_n}} N'$ and $I_s \nonfork_{N'} N_{s_n}$.
Therefore $I_s \nonfork_{M_s} N' \cup N_{s_n}$.
By Dominance $N_s \nonfork_{M_s} N' \cup N_{s_n}$.
We are done by monotonicity.
\end{proof}

We will now use DOP to construct systems as in the claim.

Let the situation be as in the first claim.
Write $p(\bar{x}, \bar{b}) = p \restriction \bar{b}$.
Let $\langle \bar{a}^\alpha_1 \bar{a}^\alpha_2 \alpha < \mu \rangle$
be a Morley sequence for $\ftp(\bar{a}_1 \bar{a}_2/M)$.
Such a Morley sequence exists by assumption on $\mathfrak{C}$ and
stationarity over models.
Let $M_i^\alpha$ be prime over $M \cup \bar{a}_i^\alpha$, for $i=1,2$.
Such a prime model exists by the axioms. 
Then $M_1^\alpha \nonfork_M M_2^\beta$ for every $\alpha < \beta$,
by Dominance.
By axiom on prime there exists $M^{\alpha \beta}$ prime over 
$M^\alpha_1\cup M_2^\beta$.
Let $\bar{b}^{\alpha \beta}$ be the image of $\bar{b}$ in $M^{\alpha \beta}$.
Let $p^{\alpha \beta} = 
p(\bar{x}, \bar{b}^{\alpha \beta})_{M^{\alpha \beta}} 
\in S(M^{\alpha \beta})$, 
which exists and is
regular since $p$ is based on $\bar{b}$.
Thus, $p^{\alpha \beta} \perp M^\alpha_1$ and 
$p^{\alpha \beta} \perp M^\beta_2$.
Let $I^{\alpha \beta}$ be a Morley sequence of length $\mu$ for
$p^{\alpha \beta}$.
Let $N^{\alpha \beta}$ be prime over $M^{\alpha \beta} \cup I^{\alpha \beta}$.
Then, for the claim, for each $R \subseteq [\mu]^2$, 
the system 
\[
\mathcal{S}_R = \{ M \} \cup
\langle M_i^\alpha : \alpha < \mu, i = 1,2 \rangle
\cup \langle N^{\alpha \beta} : \langle \alpha ,\beta\rangle \in R \rangle
\]
is an independent system.
Hence, there exists $M_R$ prime over it.

The final claim explains the name of Dimensional Order Property:
It is possible to code the relation $R$ (in particular an order in 
the following theorem) by looking at dimensions of indiscernibles
in a model $M_R$.
Note that the converse holds also, namely that the following property
characterises DOP (we do not prove this fact as it is not necessary
to obtain the main gap).
Recall $\mu > \lambda(D)$.

\begin{claim}\label{R and MR}
 The pair $\langle \alpha, \beta \rangle \in R$ if and only
if there exists $\bar{c} \in M_R$ with the property that
$\ftp(\bar{a}_1 \bar{a}_2\bar{b}/\emptyset)=
\ftp(\bar{a}^\alpha_1 \bar{a}^\beta_2\bar{c}/\emptyset)$
and for every prime $M^* \subseteq M_R$
over $M \cup \bar{a}^\alpha_1 \bar{a}^\beta_2$ containing $\bar{c}$ there
exists a Morley sequence for $p(\bar{x}, \bar{c})_{M^*}$ of length $\mu$.
\end{claim}
\begin{proof}
If the pair $\langle \alpha, \beta \rangle \in R$,
then $p^{\alpha \beta}$ is based on $\bar{b}^{\alpha \beta}$.
Furthermore, $I^{\alpha \beta}$ is a Morley sequence of length
$\mu$ for $p^{\alpha \beta}$ in $M_R$.
Let $M'$ be prime over 
$M \cup \bar{a}^\alpha_1 \bar{a}^\beta_2$ containing $\bar{b}^{\alpha \beta}$,
then $p(\bar{x}, \bar{b}^{\alpha \beta})_{M'}$ is realized by every element
of $I^{\alpha \beta}$ except possibly $\lambda(D)$ many.
Hence, there exists a Morley sequence of length $\mu$, since
$\mu > \lambda(D)$.

For the converse, 
let $\alpha < \beta < \mu$ be given such that $(\alpha, \beta) \not \in R$.
Let $t = (\alpha, \beta)$.
Let $\bar{c} \subseteq M_R$ finite as in the claim.
By using an automorphism, we have that 
$\ftp(\bar{c}/\bar{a}_1^\alpha\bar{a}_2^\beta)$ isolates 
$\ftp(\bar{c}/M_1^\alpha M_2^\beta)$ and hence
there exists $M_t \subseteq M_R$ prime over $M_1^\alpha M_2^\beta$ containing
$\bar{c}$.
By assumption on $\bar{c}$, there exists $I \subseteq M_R$ a Morley
sequence for $p(\bar{x}, \bar{c})_{M_t}$ of length $\mu$.
Let $N_t$ be prime over $M_t(I)$,
which exists by assumption on prime.
By the previous claim, the following system is independent
\[
\{ M \} \cup
\langle M_i^\alpha \mid \alpha < \mu, i = 1,2 \rangle
\cup \langle N^{\alpha \beta} \mid \langle \alpha ,\beta\rangle \in R \rangle
\cup \{N_t \}.
\]
Thus, in particular 
$N_t \nonfork_{M_t} \bigcup_{i < \mu} M_i \cup \bigcup_{s \in R} N_s$,
so 
$\bar{a} \nonfork_{M_t} \bigcup_{i < \mu} M_i \cup \bigcup_{s \in R} N_s$,
for each $\bar{a} \in I$.
By Dominance $\bar{a} \nonfork_{M_t} M_R$ and so $\bar{a} \in M_t$.
This is a contradiction. 
\end{proof}

All the technology is now in place to apply the methods
of \cite{Sh b} or \cite{GrHa} with the previous claim
and to derive:

\begin{theorem}\label{dopmany}
Suppose that $\mathcal{K}$ has DOP.
Then, $\mathcal{K}$ contains $2^\lambda$
nonisomorphic models of cardinality $\lambda$,
for each $\lambda > |D|+|T|$.
\end{theorem}

\begin{theorem}\label{dopmany.2}
Suppose that the class of $(D,\mu)$-homogeneous models
of a totally transcendental diagram $D$ has DOP.
Then, for each $\lambda > |D|+|T| +\mu$ there are $2^\lambda$
nonisomorphic $(D,\mu)$-homogeneous models of cardinality $\lambda$.
\end{theorem}

\bigskip
\section{Depth and the main gap}

We have now showed that if every model $(D,\aleph_0)$-homogeneous
model of a totally transcendental diagram $D$ has
NDOP, then every such model admits a decomposition.
We will introduce an equivalence between decompositions,
as well as the notion of depth,
in order to compute the spectrum function for $\mathcal{K}$.
Most of the treatment will be done under the assumption
that $\mathcal{K}$ has NDOP.

\begin{definition}
We say that $\mathcal{K}$ has NDOP if every $N \in \mathcal{K}$
has NDOP.
\end{definition}

We introduce the \emph{depth} of a regular type.

\begin{definition}
Let $p \in S_D(M)$ be regular.
We define the \emph{depth} of $p$, written $\D(p)$.
The depth $\D(p)$ will be an ordinal, $-1$, or
$\infinity$ and we have the usual ordering
$-1 < \alpha < \infinity$ for any ordinal $\alpha$.
We define the relation $\D(p) \geq \alpha$ 
by induction on  $\alpha$.

\begin{enumerate}

\item 
$\D(p) \geq 0$ if $p$ is regular;

\item 
$\D(p) \geq \delta$, when $\delta$ is a limit ordinal,
if $\D(p) >\alpha$ for every $\alpha < \delta$;

\item  
$\D(p) \geq \alpha + 1$  if there exists $\bar{a}$ realizing
$p$ and a regular type $r \in S_D(M(\bar{a}))$ such
that $r \perp M$ and $\D(r) \geq \alpha$.

\end{enumerate}

We write:
\item
$\D(p)=-1$ if $p$ is not regular;
\item
$\D(p)=\alpha$ if $\D(p) \geq \alpha$ but it is not
the case that $\D(p) \geq \alpha + 1$;
\item
$\D(p)=\infinity$ 
if $\D(p) \geq \alpha$ for every ordinal $\alpha$. 


We let 
$\D(\mathcal{K}) = \sup \{ \D(p)+1 \mid M \in \mathcal{K}, p \in S_D(M) \}$.
This is called the \emph{depth} of $\mathcal{K}$.
\end{definition}

\begin{lemma}
Let $p \in S_D(M)$ be regular with $\D(p) < \infinity$.
Let $\bar{a} \models p$ with $r \in S_D(M(\bar{a}))$ regular with 
$r \perp M$.
Then $\D(r) < \D(p)$.
\end{lemma}
\begin{proof} 
This is obvious, by definition of depth, if
$\D(r)=\D(p)$ is as above, then $\D(p) \geq \D(p)+1$,
contradicting $\D(p) < \infinity$.
\end{proof}

\begin{lemma}\label{initial}
Let $p \in S_D(M)$ be regular.
If $\D(p) < \infinity$ and $\alpha \leq \D(p)$,
then there
exists $q$ regular such that $\D(q) = \alpha$.
\end{lemma}
\begin{proof}
By induction on $\D(p)$. For $\D(p)=0$ it is clear.
Assume that $\D(p)= \beta+1$.
Let $\bar{a} \models p$ and
let $r \in S_D(M(\bar{a}))$ be such that 
$r \perp M$ and $\D(r) \geq \beta$.
Then, by the previous lemma, 
$\D(r)=\beta$. Hence, we are done by induction.
Assume that $\D(p)= \delta$, where $\delta$ is a limit
ordinal.
Let $\alpha < \delta$. 
Then, $\D(p) > \alpha$ by definition, so 
there exist $\bar{a} \models p$ and
$r \in S_D(M(\bar{a}))$ regular such that 
$r \perp M$ and $\D(r) \geq \alpha$. 
By the previous lemma $\D(r) < \D(p)$,
so we are done by induction.
\end{proof}

We first show that 
the depth respects the equivalence relation $\not \perp$.
\begin{lemma}
Let $p, q \in S_D(M)$ be regular such that $p \not \perp q$.
Then $\D(p)=\D(q)$.
\end{lemma}
\begin{proof}
By symmetry, it is enough to show that $\D(p)\leq \D(q)$.
We show by induction on $\alpha$ that $\D(p)\geq \alpha$
implies $\D(q) \geq \alpha$.
For $\alpha=0$ or $\alpha$ a limit ordinal, it is obvious.
Suppose that $\D(p)\geq \alpha+1$, and let $\bar{a}$ realize $p$
and $r \in S_D(M(\bar{a}))$ be such that $\D(r) \geq \alpha$
and $r \perp M$.
Since $p \not \perp q$, by Lemma \ref{strongequiv}, there
exists $\bar{b}$ realizing $q$ such that $M(\bar{a}) = M(\bar{b})$.
This implies that $\D(q) \geq \alpha+1$.
\end{proof}

\begin{lemma} \label{pardep}
Suppose $\mathcal{K}$ has NDOP.
Let $M \subseteq N$, with $M,N \in \mathcal{K}$.
Let $p \in S_D(M)$ be regular.
Then $\D(p) = \D(p_N)$.
\end{lemma}
\begin{proof}
We first show that $\D(p) \geq \D(p_N)$.
By induction on $\alpha$, we show that $\D(p) \geq \alpha$
implies $\D(p_N) \geq \alpha$.

For $\alpha=0$ it follows from the fact that $p_N$ is regular.
For $\alpha$ a limit ordinal it follows by induction.
Suppose $\D(p) \geq \alpha+1$.
Let $\bar{a}$ realize $p$
and $r \in S_D(M(\bar{a}))$ regular 
be such that $\D(r) \geq \alpha$
and $r \perp M$. 
Without loss of generality, we may assume that $\bar{a} \nonfork_M N$.
Hence, by Dominance $M(\bar{a}) \nonfork_M N$.
Since $r \perp M$, then Lemma \ref{reguperp} implies that $r \perp N$.
By induction hypothesis $\D(r_{N(\bar{a})}) \geq \alpha$.
Hence $\D(p_N) \geq \alpha+1$.

The converse uses NDOP.
We show by induction on $\alpha$ that 
$\D(p_N) \geq \alpha$ implies $\D(p)\geq \alpha$.
For $\alpha =0$ or $\alpha$ a limit ordinal, this is clear.
Suppose $\D(p_N) \geq \alpha+1$.
Let $\bar{a}$ realize $p_N$.
Then $\bar{a} \nonfork_M N$, so by Dominance $M(\bar{a}) \nonfork_M N$.
Consider $N'$ $D^s_{\aleph_0}$-primary over $M(\bar{a}) \cup N$.
We may assume that $N' = N(\bar{a})$.
Hence, there is $r \in S_D(N')$ regular 
such that $\D(r) \geq \alpha$
and $r \perp N$. 
Hence, by NDOP, we must have $r \not \perp M(\bar{a})$.
Therefore, by (Perp I) there exists a regular type
$q \in S_D(M(\bar{a})$ such that $r \not \perp q$.
But, since $r \perp M$, also $q \perp M$.
Moreover, by Parallelism, $r \not \perp q_{N'}$ and since $q_{N'}$
is regular, the previous lemma shows that $\D(q_{N'}) = \D(r) \geq \alpha$.
Hence, by induction hypothesis, $\D(q) \geq \alpha$. 
This implies that $\D(p) \geq \alpha+1$.
\end{proof}

Let $\lambda(D)= |D|+|T|$.
As we saw in Chapter III, if $D$ is totally transcendental,
then $D$ is stable in $\lambda(D)$.

\begin{lemma} \label{deeplambda}
Let $\mathcal{K}$ have NDOP.
If $\D(\mathcal{K}) \geq \lambda(D)^+$ then $\D(\mathcal{K})= \infinity$.
\end{lemma}
\begin{proof}
Let $p$ be regular based on $B$.
Let $M$ be $D^s_{\aleph_0}$-primary over the empty set.
Then $\| M \| \leq \lambda(D)$.
By an automorphism, we may assume that $B \subseteq M$.
Then, by Lemma \ref{pardep},
we have $\D(p)= \D(p \restriction M)$.
Thus, since $|S_D(M)| \leq \lambda(D)$,
there are at most $\lambda(D)$ possible depths.
By Lemma \ref{initial}, they form an initial segment
of the ordinals. 
This proves the lemma.
\end{proof}

\begin{definition}
The class $\mathcal{K}$ is called 
\emph{deep} if $\D(\mathcal{K}) = \infinity$.
\end{definition}

The next theorem is the main characterization of deep $\mathcal{K}$.
A class $\mathcal{K}$ is deep if and only if a natural partial order 
on $\mathcal{K}$ is not well-founded.
This will be used to construct nonisomorphic models in Theorem \ref{deepmany}.

\begin{theorem}\label{deep}
$\mathcal{K}$ is deep if and only if there exists a sequence
\[
\langle M_i, \bar{a}_i \mid i < \omega \rangle
\] 
such that
\begin{enumerate}
\item 
$M_0$ has cardinality $\lambda(D)$;
\item
$\ftp(\bar{a}_i/M_i)$ is regular;
\item
$M_{i+1}$ is prime over $M_i \cup \bar{a}_i$;
\item
$M_{i+1}/M_i \perp M_{i-1}$, if $i > 0$.
\end{enumerate}
\end{theorem}
\begin{proof}
Suppose that $\mathcal{K}$ is deep.
Prove by induction on $i < \omega$ that a sequence
satisfying (1)--(4) exists and that in addition
\begin{enumerate}\setcounter{enumi}{4}
\item
$\D(\ftp(\bar{a}_i/M_i))= \infinity$.
\end{enumerate}
This is possible.
For $i=0$, let $M \in \mathcal{K}$ and
$p \in S_D(M)$ be regular such that $\D(p) \geq \lambda(D)^++1$.
Such a type exists since $\mathcal{K}$ is deep.
Now, let $B$ be finite such that
$p$ is regular over $B$.
Let $M_0 \in \mathcal{K}$ contain $B$ be of cardinality $\lambda(D)$.
Then, since $p = (p\restriction M_0)_M$, we 
have $\D(p\restriction M_0)$ by Lemma \ref{pardep}.
Let $\bar{a}_0$ realize $p \restriction M_0$.
By the previous fact, $\D( \ftp(\bar{a}_0/M_0)) = \infinity$.
Now assume that $\bar{a}_i, M_i$ have been constructed.
Let $M_{i+1}$ be prime over $M_i \cup \bar{a}_i$.
By (5), we must have $\D(\ftp(\bar{a}_i/M_i)) \geq \lambda(D)^++1$,
so there exists $\bar{a}_{i+1}$ realizing $\ftp(\bar{a}_i/M_i)$
and a regular type $p_i \in S_D(M_{i+1})$ such that
$\D(p_i) \geq \lambda(D)^+$ and $p_i \perp M_i$.
Let $\bar{a}_{i+1}$ realize $p_i$,
then (1)--(5) hold.

For the converse, suppose there exists 
$\langle M_i, \bar{a}_i \mid i < \omega \rangle$ satisfying (1)--(4).
We show by induction on $\alpha$ that $\D(\ftp(\bar{a}_i/M_i)) \geq \alpha$,
for each $i < \omega$. 
This is clearly enough since then $\D(\ftp(\bar{a}_0/M_0))=\infinity$.
For $\alpha = 0$, this is given by (2), and for $\alpha$ a limit ordinal,
this is by induction hypothesis.
For the successor case, assume that $\D(\ftp(\bar{a}_i/M_i)) \geq \alpha$,
for each $i < \omega$.
Fix $i$.
Then by (4) $\ftp(\bar{a}_{i+1}/M_{i+1}) \perp M_i$.
By (2) $\ftp(\bar{a}_{i+1}/M_{i+1})$ is regular and
by (3) $M_{i+1} = M_i(\bar{a}_i)$.
By induction hypothesis $\D(\ftp(\bar{a}_{i+1}/M_{i+1})) \geq \alpha$,
hence $\D(\ftp(\bar{a}_i/M_i) \geq \alpha+1$ by definition of depth.
\end{proof}

Recall the following definition.

\begin{definition} We say that $A$ \emph{dominates $B$ over $M$} if
for every set $C$, if $A \nonfork_M C$
then $B \nonfork_M C$.
\end{definition}

We rephrase some of the results we have obtained in the following
remark.

\begin{remark}
For any set $A$, $A$ dominates $M(A)$ over $M$.
Thus, if $M \subseteq N$, and $\bar{a} \in N \setminus M$
there always is a model $M'$ such that $\bar{a} \in M' \subseteq N$
and $M'$ is maximally dominated by $\bar{a}$ over $M$,
i.e. $M'$ is dominated by $\bar{a}$ over $M$ and every
model contained in $N$ strictly containing $M'$ is \emph{not} dominated
by $\bar{a}$ over $M$.
\end{remark}

We introduce triviality.
The name comes from the fact that the pregeometry on the
set of realizations of a trivial type is trivial.

\begin{definition}
A type $p \in S_D(M)$ is \emph{trivial} if for every 
$M', N \in \mathcal{K}$ such that
$M \subseteq M' \subseteq N$ and for every set $I \subseteq p_{M'}(N)$
of pairwise independent sequences over $M'$,
then $I$ is a Morley sequence for $p_{M'}$. 

If $p$ is trivial, $\bar{a} \models p$ and
$\bar{a}$ dominates $B$ over $M$, then
we say that $B/M$ is \emph{trivial}.
\end{definition}

\begin{remark}
If $\ftp(\bar{a}/M)$ is trivial, then $M(\bar{a})/M$ is trivial.
\end{remark}

The next lemma says essentially that all the regular types
of interest are trivial.

\begin{lemma}\label{trivial} If $\mathcal{K}$ has NDOP, then
if $p \in S_D(M)$ is regular with $\D(p) > 0$, then
$p$ is trivial.
\end{lemma}

\begin{proof}
Suppose $p \in S_D(M)$ is not trivial.
Without loss of generality 
$\bar{a}_i$ for $i \leq 2$ be pairwise independent over $M$ such
that $\{ \bar{a}_i \mid i \leq 2 \}$ is not. 
Since $\D(p) > 0$, by using an automorphism, we can
find $r \in S_D(M(\bar{a}_0))$ regular such that $r \perp M$.

Let $N = M(\bar{a}_0,\bar{a}_1,\bar{a}_2)$.
Let $M' \subseteq N$ be maximal such that $\bar{a}_1 \bar{a}_2 \nonfork_M M'$.
Thus, we may assume that $N = M'(\bar{a}_0,\bar{a}_1,\bar{a}_2)$.
Since $\bar{a}_0$ realizes $p$, and $M'/M \perp p$, we have
$\bar{a}_0 \nonfork_M M'$.
Hence, by Lemma \ref{reguperp}, we must have $r \perp M'$.
By the previous remark, choose $M_i \subseteq N$ maximally dominated
by $\bar{a}_i$ over $M'$. 
By choice of $M_i$ we have $M_1 \nonfork_M M_2$.
Thus, by definition of $M_i$ and NDOP, necessarily
$N$ is $D^s_{\aleph_0}$-primary over $M_1 \cup M_2$.

Now, since $M'/M \perp p$, we have $\bar{a}_0\bar{a}_i \nonfork_M M'$.
Hence $M(\bar{a}_0) \nonfork_{M'} M_i$, for $i=1,2$.
By Lemma \ref{reguperp}, we have $r_{N} \perp M_i$ for $i=1,2$,
contradicting NDOP. 
\end{proof}

The next lemmas are used to calculate the spectrum function.

\begin{lemma} \label{bas}
Assume $\mathcal{K}$ has NDOP.
Let $\langle M_\eta \mid \eta \in J \rangle$ be a 
complete decomposition of $N^*$ over
$M$.
Let $I$ be a subtree of $J$.
Then there exists $N_I \subseteq N^*$ and
$N_\eta \subseteq N^*$ for each $\eta \in J \setminus I$
such that
$\{ N_I \} \cup \{ N_\eta \mid \eta \in J \setminus I \}$ is a complete
decomposition of $N^*$ over $N_I$.
\end{lemma}
\begin{proof}
Define $N_I \subseteq N^*$ and $N_\eta \subseteq N^*$ 
for $\eta \in J \setminus I$ as follows
\begin{enumerate}
\item
$N_I$ is $D^s_{\aleph_0}$-primary over 
$\bigcup \{M_\eta \mid \eta \in I \}$;

\item
$N_I \nonfork_{M_I} M_J$;

\item
$N_\eta = N_{\eta^-}(M_\eta)$ for $\eta \in J \setminus I$
and when $\eta^- \in I$ then $N_\eta = N_I(M_\eta)$;

\item
$N_\eta \nonfork_{M_\eta} \bigcup_{\eta \prec \nu} M_\nu$.

\end{enumerate}
This is easily done and one checks immediately that
it satisfies the conclusion of the lemma.

\end{proof}

We now define an equivalence relation on decompositions.

\begin{definition}
Let $\langle M_\eta \mid \eta \in I \rangle$ be a complete decomposition
of $N^*$.
Define an equivalence relation $\sim$ on $I \setminus \{ \langle \rangle \}$
by 
\[
\eta \sim \nu
\quad
\text{if and only if}
\quad
M_\eta / M_{\eta^-} \not\perp M_\nu / M_{\nu^-}.
\]
\end{definition}
By Equivalence, this is indeed an equivalence relation.
By the following lemma, any two sequences in the same 
$\sim$-equivalence class have a common predecessor.

\begin{lemma} \label{abov}
If $\langle M_\eta \mid \eta \in I \rangle$ is a decomposition of $N^*$,
then for $\eta, \nu \in I \setminus \{ \langle \rangle \}$
such that $\eta^- \not = \nu^-$ we have 
$M_\eta / M_{\eta^-} \perp M_\nu / M_{\nu^-}$.
\end{lemma}
\begin{proof}
Let $\eta, \nu \in I \setminus \{ \langle \rangle \}$
such that $\eta^- \not = \nu^-$. 
Let $u$ be the largest common sequence of $\eta^-$ and $\nu^-$.
We have $M_{\eta^-} \nonfork_{M_u} M_{\nu^-}$, by 
independence of the decomposition. 
By definition $M_\eta / M_{\eta^-} \perp M_{\eta^{--}}$.
Hence, by Lemma \ref{reguperp}, we have 
$M_\eta / M_{\eta^-} \perp M_u$ and also 
$M_\eta / M_{\eta^-} \perp M_{\nu^-}$.
Therefore $M_\eta / M_{\eta^-} \perp M_\nu/ M_{\nu^-}$. 
\end{proof}

The next lemma will be used inductively.

\begin{lemma}\label{induct}
Let $\langle M_\eta \mid \eta \in I \rangle$ 
and $\langle N_\nu \mid \nu \in J \rangle$ be a complete
decompositions of $N^*$ over $M$.
Let $I' = \{ \eta \in I \mid \eta^- = \langle \rangle \}$
and $J' = \{ \nu \in J \mid \nu^- = \langle \rangle \}$.
Then there exists a bijection $f \colon I' \rightarrow J'$
such that
\begin{enumerate}
\item
$f$ preserves $\sim$-classes;
\item
If $\eta \in I'$ and $M_\eta/M$ is trivial 
then $M_\eta \fork_M N_{f(\eta)}$.
\end{enumerate}
\end{lemma}
\begin{proof}
Choose a representative for each $\not \perp$-class
among the regular types of $S_D(M)$.
Build the bijection by pieces.
For each regular $p \in S_D(M)$, the cardinalities of 
$\{ \eta \in I \mid M_\eta / M \not \perp p \}$ and
$\{ \nu \in J \mid N_\nu / N \not \perp p \}$ are equal
and both equal to the dimension of $p(N^*)$ by construction.
If $p$ is not trivial, then choose any bijection between
the two sets.
If $p$ is trivial, for each $\eta \in I$ such that
$M_\eta / M \not \perp p$ there exists exactly one 
$\nu \in J'$ such that $M_\eta \fork_M N_\nu$.
Let $f$ send each such $\eta$ to their corresponding $\nu$.
Since there is no relation between $p$'s belonging
to different equivalence classes, this 
is enough.
\end{proof}

The following quasi-isomorphism will be relevant for
the isomorphism type of models.

\begin{definition}
Two $\omega$-trees $I, J$ are said to be \emph{quasi-isomorphic},
if there exists a partial function $f$ from $I$ to $J$ such that
\begin{enumerate}
\item
$f$ is order-preserving;
\item
For each $\eta \in I$ all but at most $\lambda(D)$ many
successors of $\eta$ are in $\dom(f)$;

\item
For each $\nu \in J$ all but at $\lambda(D)$ many
successors of $\nu$ are in the $\ran(f)$.
\end{enumerate}
A function $f$ as above is called a \emph{quasi-isomorphism}.
\end{definition}

\begin{theorem} \label{quasiso}
Let $\langle M_\eta \mid \eta \in I \rangle$ and
$\langle M'_\nu \mid \nu \in J \rangle$ be complete decompositions
of $N^*$.
Then there exists a $\sim$-class preserving quasi-isomorphism
from $I$ to $J$.
\end{theorem}
\begin{proof}
For each $\eta \in I$, 
let $I_\eta^+ = \{ \nu \in I \mid \nu^- = \eta \}$.
We define a partial class preserving function $f_\eta$
from $I_\eta^+$ into $J$ as follows.
Then $M_\eta$ has cardinality $\lambda(D)$, so we can find
$I_0$ and $J_0$ of cardinality at most $\lambda(D)$ such
that there exists $N \subseteq N^*$ containing $M_\eta$,
such that $M$ is $D^s_{\aleph_0}$-primary over both 
$\bigcup \{ M_\nu \mid \nu \in I_0 \}$ and 
$\bigcup \{ M'_\nu \mid \nu \in J_0 \}$.
By Lemma \ref{bas} and Lemma \ref{induct},
there exists a 
partial function $f_\eta$ from $I^+_\eta \setminus I_0$
into $J$ satisfying conditions (1) and (2) in Lemma \ref{induct}.

Now let $f = \bigcup_{\eta \in I} f_\eta$ (we let $f_{\langle \rangle}$ map
$\langle \rangle$ to $\langle \rangle$).
Clearly $f$ is well-defined, since the domains of all the
$f_\eta$'s are disjoint. 
Further, by construction, the condition involving $\lambda(D)$ is
satisfied.

It remains to show that $f$ is one-to-one and order preserving.
We check order preserving and leave one-to-one to the 
reader.
Let $\eta \prec \nu \in I$ be given.
We may assume that $\eta \not = \langle \rangle$.
Then, by Lemma \ref{trivial}, we have $M_\eta/M_{\eta^-}$ is trivial.
We are going to compute $f(\eta)$ and $f(\nu)$.
Recall that $f(\eta) = f_{\eta^-}(\eta)$.
In the notation of Lemma \ref{bas} and of the first paragraph,
we have
\[
\langle N_\zeta : \zeta \in I \setminus I_0 \rangle \cup \{ N \}
\quad
\text{and}
\quad
\langle N'_\zeta : \zeta \in J \setminus J_0 \rangle \cup \{ N \},
\]
two complete decompositions of $N^*$ over $N$.
By Lemma \ref{induct}, we have 
\[
N_{\eta^-} \fork_N N'_{f(\eta^-)}.
\]
Then, necessarily $M_{\nu}/M_{\nu^-} \not \perp M'_{f(\nu)}/M'_{f(\nu^-)}$
and any sequence $\sim$-related to $\nu$ is $\prec$-above $\eta$.
Consider the following independent tree
\[
\langle N'_\zeta :
\zeta \in I \setminus I_0, f_\eta^-(\eta) \not \prec \zeta
\rangle \cup 
\langle N_\zeta : \eta \in I, \eta \prec \zeta \rangle \cup \{ N \}.
\]
By triviality of $M_\eta/M_{\eta^-}$, it is a decomposition of $N^*$ over
$N$.
Hence, by Lemma \ref{abov} we have $M_\nu/M_{\nu^-} \perp N_\zeta/N_{\zeta^-}$,
for each $\zeta \in I \setminus I_0, f_\eta^-(\eta) \not \prec \zeta$.
This implies that the $\sim$-class of $f_{\eta^-}(\nu)$
is above $f_\nu( \eta^-)$.
Thus, $f$ is order preserving. 
\end{proof}

In order to construct many nonisomorphic models, 
we will need a special kind of trees.
For an $\omega$-tree $I$ and $\eta \in I$, 
denote by $I_\eta = \{ \nu \in I \mid \eta \prec \nu \}$.
We write $I_\eta \cong I_\nu$ if both trees are isomorphic as trees.
\begin{definition}
An $\omega$-tree $I$ is called \emph{ample} if for every $\eta \in I$,
with $\eta^- \in I$, we have 
\[
|\{ \nu \in I : \nu^- = \eta^- \text{ and } I_\nu \cong I_\eta \}| >
\lambda(D).
\]
\end{definition}

We now state a fact about ample $\omega$-trees.
If $I$ is a tree, by definition every $\eta \in I$ 
is well-founded in the order of $I$.
The \emph{rank} of $\eta$ in $I$ will be the natural rank
associated with the well-foundness relation on $\eta$ in $I$.

\begin{fact}
Let $I$, $J$ be ample trees.
Let $f$ be a quasi-isomorphism from $I$ to $J$.
Then for each $\eta \in \dom(f)$, the rank of $\eta$ in $I$
is equal to the rank of $f(\eta)$ in $J$.
\end{fact}

In the next proof, write $\ell(\eta)$ for the level of $\eta$.

\begin{theorem}\label{deepmany}
If $\mathcal{K}$ is deep, for each $\mu > \lambda(D)$,
there are $2^\mu$ 
nonisomorphic models of cardinality $\mu$.
\end{theorem}
\begin{proof}
Let $\mu > \lambda(D)$.
Since $\mathcal{K}$ is deep by Theorem \ref{deep}, there
exists 
\[
\langle M_i, \bar{a}_i \mid i < \omega \rangle
\] 
such
that
\begin{enumerate}
\item 
$M_0$ has cardinality $\lambda(D)$;
\item
$\ftp(\bar{a}_i/M_i)$ is regular;
\item
$M_{i+1}$ is prime over $M_i \cup \bar{a}_i$;
\item
$M_{i+1}/M_i \perp M_{i-1}$, if $i > 0$.
\end{enumerate}
Let $p = \ftp(\bar{a}_0/M_0)$.
Then $p$ is regular based on a finite set $B$.
We will find $2^\mu$ non-isomorphic models of size $\mu$ with
$B$ fixed.
This implies the conclusion of the theorem since $\mu^{<\aleph_0} = \mu$.

For each $X \subseteq \mu$ of size $\mu$, let $I_{X}$ be an
ample $\omega$-tree with the property that the set of ranks of
elements of the first level of $I_X$ is exactly $X$.
Such a tree clearly exists ($\mu > \lambda(D)$).
Define the following system 
$\langle M^X_\eta \mid \eta \in I_X \rangle$:
\begin{enumerate}
\item
$M^X_{\langle \rangle} = M_0$;
\item
If $\eta_0 \prec \dots \prec \eta_n \in I_X$,
then 
$$\ftp(M^X_{\eta_0} \dots M^X_{\eta_n} /\emptyset) = 
\ftp(M_{\ell(\eta_0)} \dots M_{\ell(\eta_n)}/\emptyset).$$
\end{enumerate} 
This is easy to do and by choice of 
$\langle M_i, \bar{a}_i \mid i < \omega \rangle$ this is a decomposition.
Let $M_X$ be a $D^s_{\aleph_0}$-primary model over 
$\bigcup \{ M^X_\eta \mid \eta \in I_X \}$.
Then $M_X \in \mathcal{K}$ has cardinality $\mu$.
By NDOP, $\langle M_i, \bar{a}_i \mid i < \omega \rangle$ is a complete
decomposition of $M_X$ over $M_0$.

We claim that for $X \not = Y$ as above, $M_X \not \cong_B M_Y$.
Let $X, Y \subseteq \mu$ of cardinality $\mu$ be such
that $X \not = Y$.
Suppose $M_X \cong_B M_Y$.
Then, by Theorem \ref{quasiso}, there exists a class-preserving
quasi-isomorphism
between $I_X$ and $I_Y$.
Since $B$ is fixed, the first level of $I_X$ is mapped to the
first level of $I_Y$.
By the previous fact, we conclude that $X=Y$, a contradiction.
\end{proof}

We have shown that deep diagrams have many models. 
The usual methods (see \cite{Sh b} for example) can be used
to compute the spectrum of $\mathcal{K}$ when 
$\mathcal{K}$ is \emph{not} deep.
Recall that when $\mathcal{K}$ has NDOP but is not
deep then $\D(\mathcal{K}) < \lambda(D)^+$, by Lemma \ref{deeplambda}.

\begin{theorem}\label{shallow}
If $\mathcal{K}$ has NDOP but is not deep, then for each ordinal $\alpha$
with $\aleph_{\alpha} \geq \lambda(D)$, we have
$I(\aleph_{\alpha}, \mathcal{K}) \leq 
\beth_{\D(\mathcal{K})}(\aleph_0+|\alpha|^{2^{|T|}})
< \beth_{\lambda(D)^+}(\aleph_0+|\alpha|)$.
\end{theorem}

This proves the \emph{main gap} for the class $\mathcal{K}$
of $(D,\aleph_0)$-homogeneous
models of a totally transcendental diagram $D$.

\begin{theorem}[Main Gap] \label{maingap}
Let $\mathcal{K}$ be the class of $(D,\aleph_0)$-homogeneous
models of a totally transcendental diagram $D$.
Then, either $I(\aleph_\alpha, \mathcal{K})= 2^{\aleph_\alpha}$,
for each ordinal $\alpha$ such that $\aleph_{\alpha} > |T|+|D|$,
or $I(\aleph_{\alpha}, \mathcal{K}) <
\beth_{(|T|+|D|)^+}(\aleph_0+|\alpha|)$, for each $\alpha$ such that
$\aleph_{\alpha} > |T|+|D|$.
\end{theorem}
\begin{proof}
If $\mathcal{K}$ has DOP (Theorem \ref{dopmany}) 
or has NDOP but is deep (Theorem \ref{deepmany}), then $\mathcal{K}$
has the maximum number of models. 
Otherwise, $\mathcal{K}$ has NDOP and is not deep and the bound follows
from Theorem \ref{shallow}.
\end{proof}

Similar methods using the existence of $D^s_\mu$-prime models for
totally transcendental diagrams allow us to prove the main gap for
$(D,\mu)$-homogeneous models of a totally transcendental diagram
$D$. 

\begin{theorem}
Let $\mathcal{K}$ be the class of $(D,\mu)$-homogeneous
models of a totally transcendental diagram $D$.
Then, either $I(\aleph_\alpha, \mathcal{K})= 2^{\aleph_\alpha}$,
for each ordinal $\alpha$ such that $\aleph_{\alpha} > |T|+|D|+\mu$,
or $I(\aleph_{\alpha}, \mathcal{K}) <
\beth_{(|T|+|D|)^+}(\aleph_0+|\alpha|)$, for each $\alpha$ such that
$\aleph_{\alpha} > |T|+|D|+\mu$.
\end{theorem}

Finally, similarly to \cite{GrHa} or \cite{Ha}, it is possible
to show that for $\alpha$ large enough, the function
$\alpha \mapsto I(\aleph_{\alpha}, \mathcal{K})$ is non-decreasing,
for the class $\mathcal{K}$ of $(D, \mu)$-homogeneous models
of a totally transcendental diagram $D$.

\end{document}